\documentclass[11pt]{amsart}
\usepackage{fullpage}
\usepackage[utf8]{inputenc}
\usepackage[english]{babel}
\usepackage[T1]{fontenc}
\usepackage{graphicx}
\usepackage{amsmath}
\usepackage{amsfonts}
\usepackage{graphicx}
\usepackage{amsthm}

\usepackage{color}

\newtheorem{theo}{Theorem}[section]

\newtheorem{rem}[theo]{Remark}

\newtheorem{propo}[theo]{Proposition}
\newtheorem{lemme}[theo]{Lemma}
\newtheorem{defi}[theo]{Definition}
\newtheorem{ex}[theo]{Example}

\newtheorem{hyp}[theo]{Assumptions}
\newcommand{\E}{\mathbb{E}}
\newcommand{\R}{\mathbb{R}}
\newcommand{\PP}{\mathbb{P}}

\newcommand{\N}{\mathbb{N}}

\numberwithin{equation}{section}

\author{Charles-Edouard BREHIER}

\title{Strong and weak order in averaging for SPDEs}

\address{ENS Cachan Bretagne - IRMAR, Universit\'e Rennes 1\\
Avenue Robert Schumann\\ F-35170 Bruz\\ France}

\email{charles-edouard.brehier@bretagne.ens-cachan.fr}

\date{}

\begin{document}

\begin{abstract}
We show an averaging result for a system of stochastic evolution equations of parabolic type with slow and fast time scales. We derive explicit bounds for the approximation error with respect to the small parameter defining the fast time scale. We prove that the slow component of the solution of the system converges towards the solution of the averaged equation with an order of convergence is $1/2$ in a strong sense  - approximation of trajectories -  and $1$ in a weak sense - approximation of laws. These orders turn out to be the same as for the SDE case. 
\end{abstract}

\maketitle

\section{Introduction}

In this paper, we
consider a randomly-perturbed system of reaction-diffusion equations that can be written
\begin{equation}\label{truePDE}
\begin{gathered}
\frac{\partial x^\epsilon(t,\xi)}{\partial t}=\frac{\partial^2 x^\epsilon(t,\xi)}{\partial \xi^2}+f(\xi,x^\epsilon(t,\xi),y^\epsilon(t,\xi)),\\
\frac{\partial y^\epsilon(t,\xi)}{\partial t}=\frac{1}{\epsilon}\frac{\partial^2 y^\epsilon(t,\xi)}{\partial \xi^2}+\frac{1}{\epsilon}g(\xi,x^\epsilon(t,\xi),y^\epsilon(t,\xi))+\frac{1}{\sqrt{\epsilon}}\frac{\partial \omega(t,\xi)}{\partial t},
\end{gathered}
\end{equation}
for $t\geq 0, \xi\in(0,1)$, with initial conditions $x^\epsilon(0,\xi)=x(\xi)$ and $y^\epsilon(0,\xi)=y(\xi)$, and Dirichlet boundary conditions $x^\epsilon(t,0)=x^\epsilon(t,1)=0,y^\epsilon(t,0)=y^\epsilon(t,1)=0$.
The stochastic perturbation $\frac{\partial \omega(t,\xi)}{\partial t}$ is a space-time white noise and $\epsilon>0$ is a small parameter.

Such a system presents a specific structure: while the variations of the first component {\em a priori} depend on the slow time $t$, the second component evolves with respect to the fast time $\frac{t}{\epsilon}$. These two natural time scales are coupled through the nonlinear terms in the two equations.

In this setting, the 
%
%
%
main idea of the {\em averaging principle}, see for instance \cite{Has},  is to study the behaviour of the system when $\epsilon$ tends to $0$ by exhibiting a limit equation - the so-called averaged equation -  for the slow component $x^\epsilon$, and to prove 
the  convergence of $x^\epsilon$ towards  the solution of this averaged equation. Here, we show two approximation results  - see Theorems \ref{strong} and \ref{weak} -  and give explicit order of convergence with respect to $\epsilon$.

The averaged equation comes from the asymptotic behaviour of the fast equation. Heuristically, when $t>0$ and $\epsilon\rightarrow 0$, the fast time $\frac{t}{\epsilon}$ goes to $+\infty$, so that we expect the solution of the fast equation to be quickly close to a  stochastic equilibrium (and this is  the case under the dissipativity assumptions made in this paper), and that we can replace $y^\epsilon(t,\xi)$ in the slow equation with some stationnary - in the stochastic sense - process, which leads to the definition of averaged coefficients in the slow equation.

To give precise results, it is convenient to look at the equations in an abstract setting, where system \eqref{truePDE} can be rewritten
\begin{equation}\label{eqgen}
\begin{gathered}
dX^{\epsilon}(t)=\left(AX^{\epsilon}(t)+F(X^{\epsilon}(t),Y^{\epsilon}(t))\right)dt\\
dY^{\epsilon}(t)=\frac{1}{\epsilon}\left(BY^{\epsilon}(t)+G(X^{\epsilon}(t),Y^{\epsilon}(t))\right)dt+\frac{1}{\sqrt{\epsilon}}dW(t),
\end{gathered}
\end{equation}
with initial conditions given by $X^\epsilon(0)=x\in H$, $Y^\epsilon(0)=y\in H$, where $H$ is the Hilbert space $L^2(0,1)$, and $W$ is a cylindrical Wiener process on $H$ - see Section \ref{rappelsWiener}. In the case of system \eqref{truePDE}, the definitions of $A$ and $B$ are given in Example \ref{exampleAB}, and the definitions of $F$ and $G$ as Nemytskii operators are given in the second part of Example \ref{exFG}. Nevertheless the abstract setting allows for more general equations, and in the sequel we only work with system \eqref{eqgen}.

When $\epsilon$ tends to $0$, the slow component $X^\epsilon$ is approximated by the process $\overline{X}$, which follows the deterministic evolution equation
\begin{equation}\label{eqmoy}
d\overline{X}(t)=(A\overline{X}(t)+\overline{F}(\overline{X}(t)))dt,
\end{equation}
with the initial condition $\overline{X}(0)=x$, where the nonlinear coefficient $\overline{F}$ is obtained via an averaging procedure - explained in detail in Section \ref{known}.

In this article, we analyse the error between $X^\epsilon(t)$ and $\overline{X}(t)$, with two different criterions. We focus on the order of convergence, i.e. we bound the error by $C\epsilon^\Lambda$, where $C$ is a constant and $\Lambda$ is the order, which gives an idea of the speed of convergence of $X^\epsilon(t)$ towards $\overline{X}(t)$. As a result, we control the error made when $X^\epsilon(t)$ is approximated by $\overline{X}(t)$. For instance, the order of convergence is crucial for the analysis of numerical schemes used to approximate the slow component $X^\epsilon$. In a forthcoming work, we extend a numerical scheme for SDEs analysed in \cite{E-L-V} for systems of SPDEs satisfying the same structure assumptions as system \eqref{truePDE}. This scheme - called the Heterogeneous Multiscale Method -  is deeply based on the averaging principle: instead of computing $X^\epsilon$, we approximate $\overline{X}$ - and we can control the error we make. Moreover, the nonlinear averaged coefficient $\overline{F}$ is never explicitly calculated in the scheme, but only approximated by using numerical approximations of the values of the fast component at large times. The theorems we prove here allow to analyse the convergence of such a scheme with the same kind of criterions, and without knowing the order of convergence in the averaging principle it would not be possible to control the error made in the numerical approximation.

The two main theorems give bounds on the error between $X^\epsilon$ and $\overline{X}$; they need different dissipativity conditions \eqref{hypLg} and \eqref{hypdiss} which determine how the fast equation converges to its equilibrium, as explained below.

First, when Assumption \ref{strictdiss} holds, the error can be estimated in a strong sense, where trajectories of the processes are compared at a given time $t$:

\begin{theo}[Strong-order]\label{strong}
Assume \eqref{hypLg}. For any $0<r<1/2$, $T>0$, $x\in H$, $y\in H$, there exists $C=C(T,r,x,y)>0$ - depending also on the constants of the problem - such that for any $\epsilon>0$ and $0\leq t\leq T$
\begin{equation}
\E|X^{\epsilon}(t)-\overline{X}(t)|_{H}\leq
C\epsilon^{1/2-r}.
\end{equation}
\end{theo}
The error can also be estimated in a weak sense, where we are interested in the distance between the laws of the processes at a given time $t$; then only Assumption \ref{weakdiss} is necessary, since we only need consequences of dissipativity at the level of the transition semi-group, instead of trajectories:
\begin{theo}[Weak-order]\label{weak}
Assume \eqref{hypdiss}. For any $0<r<1$, $T>0$, $0<\theta\leq 1$, $x\in D(-A)^\theta$, $y\in H$, $\phi\in\mathcal{C}_{b}^{2}(H)$, there exists $C>0$, depending on $r$, $T$, $\phi$, $|x|_{(-A)^\theta}$, $|y|$ and the constants of the problem, such that for any $\epsilon>0$ and $t\leq T$
\begin{equation}
|\E[\phi(X^\epsilon(t))]-\E[\phi(\overline{X}(t))]|\leq C\epsilon^{1-r}.
\end{equation}
\end{theo}
The domains $D(-A)^\theta$ are usually the classical Sobolev spaces $H^{2\theta}$ with respect to the eigenbasis of $A$ - see Definition \ref{defpowers}. We remark that for the first theorem no regularity is needed for the initial condition - i.e. we can take $\theta=0$ - while we require $\theta>0$ for the second one; this is explained in the proof of Theorem \ref{weak}. We need to take a small parameter $r>0$, which can be as small as possible, but different from $0$. This is an effect of the infinite dimensional setting.

As a consequence, we can say that the strong order in averaging is $1/2$, while the weak order is $1$. It is a general fact that the weak order is greater than the strong order (since test functions $\phi$ in the Theorem are Lipschitz continuous), but it is worth proving that there is a gap; this fact was known for SDEs see (\cite{Has}, \cite{HasYin}), but had not been proved yet for SPDEs.

The strong convergence Theorem \ref{strong} is proved when the fast equation satisfies a strict dissipativity assumption: for any $x\in H$, the function $G(x,.)$ is Lipschitz continuous, with constant $L_g$ - independent of $x$ - satisfying the following condition:
$$L_g<\mu,$$
where $\mu$ is the smallest eigenvalue of the linear operator $-B$.
Thanks to this assumption, we can easily analyze the asymptotic behaviour of the fast equation with frozen slow component; we can identify a unique invariant probability measure - depending on $x$ - and show some exponential convergence to equilibrium. More precisely, we control in a strong sense the difference between two solutions of this fast equation starting from different initial conditions, and driven by the same noise $W$: under the previous assumption, the ergodicity comes from properties of the deterministic equation only. 

The weak convergence Theorem \ref{weak} needs a weaker dissipativity assumption \eqref{weakdiss}, which yields the same ergodicity properties - unique invariant probability measure, exponential convergence to equilibrium - but with different arguments: the asymptotic behaviour of the transition semi-group can be analyzed, thanks to the non-degeneracy of the noise - leading to a Strong Feller Property. A coupling method - adapted from the study of Markov processes, like in \cite{ku-shi} or \cite{matt} - implies that the laws - instead of trajectories - of the fast process issued from two different initial conditions are exponentially closed. We refer to Section \ref{known} for a precise result, and to \cite{deb-hu-tess} for a detailed proof. It seems that for the first time an averaging result is obtained for SPDEs under a weak dissipativity condition.


Notice that we have assumed that the slow equation has no white noise term $dw(t)$; as a consequence, the averaged equation is a deterministic parabolic partial differential equation. Considering a more general situation with some additive noise terms in the slow equation, independent of the noise in the fast equation, we could still prove in a similar way a strong order result, the only changes being time regularity of solutions. We would obtain order $1/5$, which can also be compared with the order $1/3$ obtained for SDEs. But if we introduce noise in the slow equation, the method we used to prove the weak order theorem becomes more complicated, and we have not extended the result to this situation so far.

In the case of stochastic differential equations, averaging results are already well-known - see for instance \cite{FreiWen}, \cite{Has}. Convergence in law or in probability of $X^\epsilon$ to $\overline{X}$ in the space $\mathcal{C}([0,T],H)$ can be shown by different techniques: by using a Hasminskii technique based on a subdivision of the interval $[0,T]$ (see \cite{E-L-V}, \cite{Liu}); a Poisson equation (see \cite{PavStu}); the method of perturbed test functions and of a martingale problem approach (see \cite{FouGarPapSol}); or an asymptotic expansion of the solutions of Kolmogorov equations (see \cite{E-L-V}, \cite{HasYin}).

As far as stochastic partial differential equations are concerned, in \cite{Ce} both the Hasminskii technique and a martingale problem approach are used; in \cite{Ce-F} a modified Poisson equation is the essential tool. Then convergence in law or in probability of $X^\epsilon$ to $\overline{X}$ in the space $\mathcal{C}([0,T],H)$ (the space of continuous functions from $[0,T]$ to $H$) is proved; but order of convergence was never given.

Our proof of Theorem \ref{strong} relies on the Hasminskii technique already known for SDEs: we introduce an auxiliary process for which the slow component of the fast variable is frozen on small intervals of a subdivision. We use H\"older regularity of order $1-r$ in time of the slow component, for which we do not need $\theta>0$. 

To prove Theorem \ref{weak}, we adapt the method of finding an expansion with respect to $\epsilon$ of the solutions of the Kolmogorov equations related to our system. This seems to be the first time that such a method is used to prove an averaging result for SPDEs. New technical difficulties due to infinite dimension arise: we use non bounded linear operators and non smooth nonlinear coefficients, and the Kolmogorov equations are more difficult to use. For these reasons, we use a reduction to finite dimension technique, keeping in mind that bounds must be independent of dimension, so that precise estimates are needed for each term appearing in the expansion. We interpret the necessity of $\theta>0$ with a singularity which needs to be integrable.

In Section \ref{prel}, we set the notations and give some results on the fast equation, allowing to define the averaged equation; we also precise the assumptions needed to prove the Theorems. Then in Section \ref{strongproofsect}, we prove the strong-order result. In Section \ref{weakproofsect}, we give the details of the method for proving the weak-order result. Finally in Section \ref{proffestimsect} and in the Appendix, we prove all the necessary estimates.
\section{Preliminaries}\label{prel}

\subsection{Assumptions and notations}

\subsubsection{Test functions}
To study weak convergence, we use test functions $\phi$ in the space $\mathcal{C}_{b}^{2}(H,\R)$ of functions from $H$ to $\R$ that are twice continuously differentiable, with first and second order bounded derivatives.

In the sequel, we often identify the first derivative $D\phi(x)\in \mathcal{L}(H,\R)$ with the gradient in $H$, and the second derivative $D^2\phi(x)$ with a linear operator on $H$, via the formulas:
\begin{gather*}
<D\phi(x),h>=D\phi(x).h \text{ for every }h\in H\\
<D^2\phi(x).h,k>=D^2\phi(x).(h,k) \text{ for every }h,k\in H.
\end{gather*}

\subsubsection{Stochastic integration in Hilbert spaces}\label{rappelsWiener}

In this section, we recall the definition of the cylindrical Wiener process and of stochastic integral on a separable Hilbert space $H$ (its norm is denoted by $|.|_{H}$ or just $|.|$). For more details, see \cite{DaP-Z1}.

We first fix a filtered probability space $(\Omega,\mathcal{F},(\mathcal{F}_t)_{t\geq 0},\PP)$. A cylindrical Wiener process on $H$ is defined with two elements:
\begin{itemize}
\item a complete orthonormal system of $H$, denoted by  $(q_i)_{i\in I}$, where $I$ is a subset of $\N$;
\item a family $(\beta_i)_{i\in I}$ of independent real Wiener processes with respect to the filtration $((\mathcal{F}_t)_{t\geq 0})$:
\end{itemize}

\begin{equation}\label{defWiener}
W(t)=\sum_{i\in I}\beta_i(t)q_i.
\end{equation}

When $I$ is a finite set, we recover the usual definition of Wiener processes in the finite dimensional space $\R^{|I|}$. However the subject here is the study of some Stochastic Partial Differential Equations, so that in the sequel the underlying Hilbert space $H$ is infinite dimensional; for instance when $H=L^{2}(0,1)$, an example of complete orthonormal system is $(q_k)=(\sin(k.))_{k\geq 1}$ - see Example \ref{exampleAB}.

A fundamental remark is that the series in \eqref{defWiener} does not converge in $H$; but if a linear operator $\Psi:H\rightarrow K$ is Hilbert-Schmidt, then $\Psi W(t)$ converges in $L^2(\Omega,H)$ for any $t\geq 0$.

We recall that a linear operator $\Psi:H\rightarrow K$ is said to be Hilbert-Schmidt when
$$|\Psi|_{\mathcal{L}_{2}(H,K)}^{2}:=\sum_{k=0}^{+\infty}|\Psi(q_k)|_{K}^{2}<+\infty,$$
where the definition is independent of the choice of the orthonormal basis $(q_k)$ of $H$.
The space of Hilbert-Schmidt operators from $H$ to $K$ is denoted $\mathcal{L}_{2}(H,K)$; endowed with the norm $|.|_{\mathcal{L}_{2}(H,K)}$ it is an Hilbert space.

The stochastic integral $\int_{0}^{t}\Psi(s)dW(s)$ is defined in $K$ for predictible processes $\Psi$ with values in $\mathcal{L}_2(H,K)$ such that $\int_{0}^{t}|\Psi(s)|_{\mathcal{L}_2(H,K)}^{2}ds<+\infty$ a.s; moreover when $\Psi\in L^2(\Omega\times[0,t];\mathcal{L}_2(H,K))$, the following two properties hold:
\begin{gather*}
\E|\int_{0}^{t}\Psi(s)dW(s)|_{K}^{2}=\E\int_{0}^{t}|\Psi(s)|_{\mathcal{L}_2(H,K)}^{2}ds, \text{ (It\^o isometry),}\\
\E\int_{0}^{t}\Psi(s)dW(s)=0.
\end{gather*}
A generalization of It\^o formula also holds - see \cite{DaP-Z1}.

For instance, if $v=\sum_{k\in\N}v_kq_k\in H$, we can define $$<W(t),v>=\int_{0}^{t}<v,dW(s)>=\sum_{k\in\N}\beta_k(t)v_k;$$
we then have the following space-time white noise property
$$\E<W(t),v_1><W(s),v_2>=t\wedge s<v_1,v_2>.$$

Therefore to be able to integrate a process with respect to $W$ requires some strong properties on the integrand; in our SPDE setting, the Hilbert-Schmidt properties follow from the assumptions made on the linear coefficients of the equations.

\subsubsection{Assumptions on the linear operators}

We have to specify some properties of the linear operators $A$ and $B$ coming into the definition of system \eqref{eqgen}; we assume that the linear parts are of parabolic type, with space variable $\xi\in(0,1)$.

We assume that $A$ and $B$ are unbounded linear operators, with domains $D(A)$ and $D(B)$, which satisfy the following assumptions:

\begin{hyp}\label{hypAB}
\begin{enumerate}
\item We assume that $(e_k)_{k\in \N}$ and $(f_k)_{k\in \N}$ are orthonormal basis of $H$, and $(\lambda_k)_{k\in \N}$ and $(\mu_k)_{k\in \N}$ are non-decreasing sequences of real positive numbers such that:
\begin{gather*}
Ae_{k}=-\lambda_{k}e_{k}\text{ for all } k\in\N\\
Bf_{k}=-\mu_{k}f_{k}\text{ for all } k\in\N.
\end{gather*}
We use the notations $\lambda:=\lambda_0>0$ and $\mu:=\mu_0>0$ for the smallest eigenvalues of $A$ and $B$.
\item The sequences $(\lambda_k)$ and $(\mu_k)$ go to $+\infty$; moreover we have some control of the behaviour of $(\mu_k)$ given by:
$$
\sum_{k=0}^{+\infty}\frac{1}{\mu_{k}^{\alpha}}<+\infty\Leftrightarrow \alpha>1/2.$$
\end{enumerate}
\end{hyp}

\begin{ex}\label{exampleAB}
$A=B=\frac{d^2}{dx^2}$, with domain $H^2(0,1)\cap H_{0}^{1}(0,1)\in L^2(0,1)$ - homogeneous Dirichlet boundary conditions: in that case $\lambda_k=\mu_k=\pi^2 k^2$, and $e_k(\xi)=f_k(\xi)=\sqrt{2}\sin(k\pi\xi)$ - see \cite{Brezis}.
\end{ex}


In the abstract setting, powers of $-A$ and $-B$, with their domains can be easily defined:

\begin{defi}\label{defpowers}
For $a,b\in[0,1]$, we define the operators $(-A)^a$ and $(-B)^b$ by
\begin{gather*}
(-A)^a x=\sum_{k=0}^{\infty}\lambda_{k}^{a}x_ke_k\in H,\\
(-B)^b y=\sum_{k=0}^{\infty}\mu_{k}^{b}y_kf_k\in H,
\end{gather*}

with domains
\begin{gather*}
D(-A)^a=\left\{x=\sum_{k=0}^{+\infty}x_ke_{k}\in H; |x|_{(-A)^a}^{2}:=\sum_{k=0}^{+\infty}(\lambda_k)^{2a}|x_k|^2<+\infty\right\};\\
D(-B)^b=\left\{y=\sum_{k=0}^{+\infty}y_kf_{k}\in H, |y|_{(-B)^b}^{2}:=\sum_{k=0}^{+\infty}(\mu_k)^{2b}|y_k|^2<+\infty\right\}.
\end{gather*}
\end{defi}

The domains $D(-A)^a$ are related to Sobolev spaces $H^{2a}(0,1)$: therefore when $x$ belongs to a space $D(-A)^a$, the exponent $a$ represents some regularity of the function $x$. 



The semi-groups $(e^{tA})_{t\geq 0}$ and $(e^{tB})_{t\geq 0}$ can be defined by the Hille-Yosida Theorem (see \cite{Brezis}). We use the following spectral formulas: if $x=\sum_{k=0}^{+\infty}x_ke_k\in H$ and $y=\sum_{k=0}^{+\infty}y_kf_k\in H$, then for any $t\geq 0$
$$e^{tA}x=\sum_{k=0}^{+\infty}e^{-\lambda_k t}x_ke_k \quad \text{and} \quad e^{tB}y=\sum_{k=0}^{+\infty}e^{-\mu_k t}y_kf_k.$$

For any $t\geq 0$, $e^{tA}$ and $e^{tB}$ are continuous linear operators in $H$, with respective operator norms $e^{-\lambda t}$ and $e^{-\mu t}$. The semi-group $(e^{tA})$ is used to define the solution $Z(t)=e^{tA}z$ of the linear Cauchy problem
$$\frac{dZ(t)}{dt}=AZ(t)\quad \text{with} \quad Z(0)=z.$$

To define solutions of more general PDEs of parabolic type, we use mild formulation, and Duhamel principle.

These semi-groups enjoy some smoothing properties that we often use in this work. Basically we need the following properties, which are easily proved using the above spectral properties. We write them for $A$, but they also hold with $B$.
\begin{propo}\label{proporegul}
Under Assumption \ref{hypAB}, for any $\sigma\in[0,1]$, there exists $C_\sigma>0$ such that we have:
\begin{enumerate}
\item for any $t>0$ and $x\in H$
$$|e^{tA}x|_{(-A)^\sigma}\leq C_\sigma t^{-\sigma}e^{-\frac{\lambda}{2}t}|x|_{H}.$$
\item for any $0<s<t$ and $x\in H$
$$|e^{tA}x-e^{sA}x|_{H}\leq C_\sigma\frac{(t-s)^\sigma}{s^\sigma}e^{-\frac{\lambda}{2}s}|x|_H.$$
\item for any $0<s<t$ and $x\in D(-A)^\sigma$
$$|e^{tA}x-e^{sA}x|_{H}\leq C_\sigma(t-s)^\sigma e^{-\frac{\lambda}{2}s}|x|_{(-A)^\sigma}.$$
\end{enumerate}
\end{propo}

Under the previous assumptions on the linear coefficients, it is easy to show that the following stochastic integral is well-defined in $H$, for any $t\geq 0$:
\begin{equation}\label{stoconvWB}
W^B(t)=\int_{0}^{t}e^{(t-s)B}dW(s).
\end{equation}

It is called a stochastic convolution, and it is the unique mild solution of
$$dZ(t)=BZ(t)dt+dW(t)\quad \text{with} \quad Z(0)=0.$$

Under the second condition of Assumption \ref{hypAB}, there exists $\delta>0$ such that for any $t>0$ we have $\int_{0}^{t}\frac{1}{s^\delta}|e^{sB}|_{\mathcal{L}_{2}(H)}^{2}ds<+\infty$; it can then be proved that $W^B$ has continuous trajectories - via the \textit{factorization method}, see \cite{DaP-Z1} - and that for any $1\leq p<+\infty$ $\sup_{t\geq 0}\E|W^B(t)|_{H}^{p}<+\infty$.

\subsubsection{Assumptions on the nonlinear coefficients}
We now give the Assumptions on the nonlinear coefficients $F,G:H\times H\rightarrow H$. First, we need some regularity properties:
\begin{hyp}\label{hypF}
We assume that there exists $0\leq \eta<\frac{1}{2}$ and a constant $C$ such that the following directional derivatives are well-defined and controlled:
\begin{itemize}
\item For any $x,y\in H$ and $h\in H$, $|D_xF(x,y).h|\leq C|h|_{H}$ and $|D_yF(x,y).h|\leq C|h|_{H}$.
\item For any $x,y\in H$, $h\in H$, $k\in D(-A)^{\eta}$, $|D_{xx}^{2}F(x,y).(h,k)|\leq C|h|_{H}|k|_{(-A)^{\eta}}$.
\item For any $x,y\in H$, $h\in H$, $k\in D(-B)^{\eta}$, $|D_{yy}^{2}F(x,y).(h,k)|\leq C|h|_{H}|k|_{(-B)^{\eta}}$.
\item For any $x,y\in H$, $h\in H$, $k\in D(-B)^{\eta}$, $|D_{xy}^{2}F(x,y).(h,k)|\leq C|h|_{H}|k|_{(-B)^{\eta}}$.
\item For any $x,y\in H$, $h\in D(-A)^{\eta}$, $k\in H$, $|D_{xy}^{2}F(x,y).(h,k)|\leq C|h|_{(-A)^{\eta}}|k|_{H}$.
\end{itemize}
We moreover assume that $F$ is bounded.
\end{hyp}

\begin{rem}
We warn the reader that constants may vary from line to line during the proofs, and that in order to use lighter notations we usually forget to mention dependence on the parameters. We use the generic notation $C$ for such constants.
\end{rem}

We assume that the fast equation is a gradient system: for any $x$ the nonlinear coefficient $G(x,.)$ is the derivative of some potential $U$. We also assume regularity assumptions as for $F$.
\begin{hyp}\label{hypG}
The function $G$ is defined through $G(x,y)=\nabla_yU(x,y)$, for some potential $U:H\times H\rightarrow \R$. Moreover we assume that $G$ is bounded, and that the regularity assumptions given in the Assumption \ref{hypF} are also satisfied for $G$.
\end{hyp}

Finally, we need to assume some dissipativity of the fast equation. Assumption \ref{strictdiss} is necessary to prove Theorem \ref{strong}, while Assumption \ref{weakdiss} is weaker and is sufficient to prove Theorem \ref{weak}.

\begin{hyp}[Strict dissipativity]\label{strictdiss}
Let $L_g$ denote the Lipschitz constant of $G$ with respect to its second variable; then
\begin{equation}\label{hypLg}\tag{SD}
L_{g}<\mu,
\end{equation}
where $\mu$ is defined in Assumption \ref{hypAB}.
\end{hyp}

\begin{hyp}[Weak Dissipativity]\label{weakdiss}
There exist $c>0$ and $C>0$ such that for any $x\in H$ and $y\in D(B)$
\begin{equation}\label{hypdiss}\tag{WD}
<By+G(x,y),y>\leq -c|y|^2+C.
\end{equation}
\end{hyp}


Indeed, according to Assumptions \ref{hypAB} and \ref{hypG}, the weak dissipativity Assumption \ref{weakdiss} is always satisfied, while strict dissipativity requires a condition on the Lipschitz constant of $G$.

\begin{ex}\label{exFG}
We give some fundamental examples of nonlinearities for which the previous assumptions are satisfied:
\begin{itemize}
\item Functions $F,G:H\times H\rightarrow H$ of class $\mathcal{C}^2$, bounded and with bounded derivatives, such that $G(x,y)=\nabla_yU(x,y)$ and satisfying \eqref{hypLg} fit in the framework, with the choice $\eta=0$.
\item Functions $F$ and $G$ can be \textbf{Nemytskii} operators: let $f:(0,1)\times \R^2\rightarrow \R$ be a measurable function such that for almost every $\xi\in(0,1)$ $f(\xi,.)$ is twice continuously differentiable, bounded and with uniformly bounded derivatives.
Then $F$ is defined for every $x,y\in H=L^2(0,1)$ by
$$F(x,y)(\xi)=f(\xi,x(\xi),y(\xi)).$$
For $G$, we assume that there exists a function $g$ with the same properties as $f$ above, such that $G(x,y)(\xi)=g(\xi,x(\xi),y(\xi))$. The strict dissipativity Assumption \ref{hypLg} is then satisfied when $$\sup_{\xi\in(a,b),x\in \R,y\in \R}|\frac{\partial g}{\partial y}(\xi,x,y)|<\mu.$$
The conditions in Assumption \ref{hypF} are then satisfied for $F$ and $G$ as soon as there exists $\eta<1/2$ such that $D(-A)^\eta$ and $D(-B)^\eta$ are continuously embedded into $L^\infty(0,1)$ - it is the case for $A$ and $B$ given in Example \ref{exampleAB}, with $\eta>1/4$.
\end{itemize}
\end{ex}

We remark that under Assumption \ref{hypG}, if we define $U_0(x,y)=\int_{0}^{1}<G(x,sy),y>ds$, we have $U_0(x,y)=U(x,y)-U(x,0)$; therefore $U_0$ is another potential for $G$, and is the only one such that for any $x\in H$ we have $U_0(x,0)=0$. In the sequel, it is therefore not restrictive to assume $U=U_0$.

Now we can define solutions of system \eqref{eqgen}; under Assumptions \ref{hypAB}, \ref{hypF}, \ref{hypG}, we notice that the nonlinearities $F$ and $G$ are Lipschitz continuous, and the following Proposition is classical - see \cite{DaP-Z1}:
\begin{propo}
For every $\epsilon>0$, $T>0$, $x\in H$, $y\in H$, system \eqref{eqgen} admits a unique mild solution $(X^{\epsilon},Y^{\epsilon})\in (\text{L}^2(\Omega,\mathcal{C}([0,T],H)))^2$:
\begin{equation}\label{eqmild}
\begin{gathered}
X^{\epsilon}(t)=e^{tA}x+\int_{0}^{t}e^{(t-s)A}F(X^{\epsilon}(s),Y^{\epsilon}(s))ds\\
Y^{\epsilon}(t)=e^{\frac{t}{\epsilon}B}y+\frac{1}{\epsilon}\int_{0}^{t}e^{\frac{(t-s)}{\epsilon}B}G(X^{\epsilon}(s),Y^{\epsilon}(s))ds
+\frac{1}{\sqrt{\epsilon}}\int_{0}^{t}e^{\frac{(t-s)}{\epsilon}B}dW(s).
\end{gathered}
\end{equation}
\end{propo}

In other words, system \eqref{eqgen} is well-posed for any $\epsilon>0$, on any finite time interval $[0,T]$.

Some properties - bounds on moments, space and time regularity, differentiability with respect to the parameters - of $X^\epsilon$ and $Y^\epsilon$ are given in the Appendix.

\subsection{Known results about the fast equation and the averaged equation}\label{known}

In this section, we just recall without proof the main results on the fast equation with frozen slow component and on the averaged equation, defined below. Proofs can be found in \cite{Ce-F} for the strict dissipative case, and the extension to the weakly dissipative situation relies on arguments explained below.

If $x\in H$, we define an equation on the fast variable where the slow variable is fixed and equal to $x$:
\begin{equation}\label{eqfig}
\begin{gathered}
dY_x(t,y)=(BY_x(t,y)+G(x,Y_x(t,y)))dt+dW(t),\\
Y_x(0,y)=y.
\end{gathered}
\end{equation}
This equation admits a unique mild solution, defined on $[0,+\infty[$.

Since $Y^\epsilon$ is involved at time $t>0$, heuristically we need to analyse the properties of $Y_x(\frac{t}{\epsilon},y)$, with $\epsilon\rightarrow 0$, and by a change of time we need to understand the asymptotic behaviour of $Y_x(.,y)$ when time goes to infinity.

Under the strict dissipativity Assumption \ref{strictdiss}, we obtain a contractivity of trajectories issued from different initial conditions and driven by the same noise:
\begin{propo}\label{strictconvexp}
With \eqref{hypLg}, for any $t\geq0$, $x,y_1,y_2\in H$  we have
$$|Y_{x}(t,y_1)-Y_{x}(t,y_2)|_H\leq e^{-\frac{(\mu-L_g)}{2}t}|y_1-y_2|_H.$$
\end{propo}
Under the weak dissipativity Assumption \ref{weakdiss}, we obtain such an exponential convergence result for the laws instead of trajectories. The proof of this result is not staightforward, and can be found in \cite{deb-hu-tess}.
\begin{propo}\label{propoexpy1y2}
With \eqref{hypdiss}, there exist $c>0$, $C>0$ such that for any bounded test function $\phi$, any $t\geq 0$ and any $y_1,y_2\in H$
\begin{equation}\label{cvexpy1y2}
|\E\phi(Y(t,y_1))-\E\phi(Y(t,y_2))|\leq C\|\phi\|_{\infty}(1+|y_1|^2+|y_2|^2)e^{-ct}.
\end{equation}
\end{propo}

The idea of coupling relies on the following formula: if $\nu_1$ and $\nu_2$ are two probability measures on a state space $S$, their total variation distance satisfies
$$d_{TV}(\nu_1,\nu_2)=\inf\left\{\PP(X_1\neq X_2)\right\},$$
which is an infimum over random variables $(X_1,X_2)$ defined on a same probability space, and such that $X_1\sim\nu_1$ and $X_2\sim \nu_2$.

The principle is to define a coupling $(Z_1(t,y_1,y_2),Z_2(t,y_1,y_2))_{t\geq 0}$ for the processes $(Y(t,y_1)_{t\geq 0}$ and $Y((t,y_2))_{t\geq 0}$ such that the coupling time $\mathcal{T}$ of $Z_1$ and $Z_2$ - i.e. the first time the processes are equal - has an exponentially decreasing tail.

This technique was first used in the study of the asymptotic behaviour of Markov chains - see \cite{Brem}, \cite{Doe}, \cite{Lin}, \cite{Me-Twee} - and was later adapted for SDEs and more recently for SPDEs - see for instance \cite{ku-shi}, \cite{matt}.

As a consequence, we can show that there exists a unique invariant probability measure associated with $Y_x$, and that the convergence to equilibrium is exponentially fast.

First, let $\nu=\mathcal{N}(0,(-B)^{-1}/2)$ be the centered Gaussian probability measure on $H$ with the covariance operator $(-B)^{-1}/2$ - which is positive and trace-class, thanks to Assumption \ref{hypAB}.

Then $\mu^x$ defined by
\begin{equation}\label{definvmeas}
\mu^x(dy)=\frac{1}{Z(x)}e^{2U(x,y)}\nu(dy),
\end{equation}
where $Z(x)\in]0,+\infty[$ is a normalization constant, is the unique probability invariant measure associated to $Y_x$.
This expression comes from the gradient structure of equation \eqref{eqfig}, given in Assumption \ref{hypG}.

Second, under both dissipativity assumptions, the convergence to equilibrium is exponential in the following sense: 
\begin{propo}\label{convexp}
If we assume \eqref{hypLg} or \eqref{hypdiss}, there exist constants $C,c>0$ such that for any bounded function $\phi:H\rightarrow \R$ or $\phi:H\rightarrow H$, $t\geq 0$ and $x,y\in H$ we have
$$|\E\phi(Y_x(t,y))-\int_{H}\phi(z)\mu^x(dz)|\leq C\|\phi\|_{\infty}(1+|y|_{H}^{2})e^{-ct}.$$
\end{propo}

Under the strict dissipativity Assumption \ref{strictdiss}, this is a consequence of Proposition \ref{strictconvexp} - see Theorem $3.5$ and Remark $3.6$ of \cite{Ce-F}; under the weak dissipativity Assumption \ref{weakdiss}, Proposition \ref{convexp} is a consequence of Proposition \ref{propoexpy1y2} and of the properties of the invariant measures $\mu^x$ - which have finite moments of any order, uniformly bounded with respect to $x$.


Now we define the averaged equation. First we define the averaged nonlinear coefficient $\overline{F}$:
\begin{defi}For any $x\in H$,
\begin{equation}\label{deffbar}
\overline{F}(x)=\int_{H}F(x,y)\mu^x(dy).
\end{equation}
\end{defi}

Using Assumptions \ref{hypF}, \ref{hypG} and the expression of $\mu^x$, we can easily prove the following properties on $\overline{F}$:
%
\begin{propo}\label{FbarLip}
There exists $0\leq \eta<1$ and a constant $C$ such that the following directional derivatives of $\overline{F}$ are well-defined and controlled:
\begin{itemize}
\item For any $x\in H$, $h\in H$, $|D\overline{F}(x).h|\leq C|h|_{H}$.
\item For any $x\in H$, $h\in H$, $k\in D(-A)^\eta$, $|D^2\overline{F}(x).(h,k)|\leq C|h|_{H}|k|_{(-A)^\eta}$.
\end{itemize}
Moreover, $\overline{F}$ is bounded and Lipschitz continuous.
\end{propo}

\begin{rem}
Even when $F$ and $G$ are Nemytskii operators, $\overline{F}$ is not such an operator in general.
\end{rem}
Then the averaged equation - see \eqref{eqmoy} in the introduction - can be defined:
$$d\overline{X}(t)=(A\overline{X}(t)+\overline{F}(\overline{X}(t)))dt,$$
with initial condition $\overline{X}(0)=x\in H$.
For any $T>0$, this deterministic equation admits a unique mild solution $\overline{X}\in\mathcal{C}([0,T],H)$.

\section{Proof of the strong-order result}\label{strongproofsect}

The main idea - inspired by the work on SDEs of Khasminskii in \cite{Has} - of the proof of Theorem \ref{strong} is the construction of auxiliary processes $(\tilde{X}^\epsilon,\tilde{Y}^\epsilon)$ for any $\epsilon$, for which the analysis is simpler.

In this section, we assume that dissipativity of the fast equation is strict: we have \ref{hypLg}.

Let $T>0$, $x\in H$, $y\in H$ and $\epsilon>0$ be fixed. We introduce the parameter
\begin{equation}\label{choiceepsilon}
\delta=\delta(\epsilon)=\sqrt{\epsilon}
\end{equation}
to define a subdivision of $[0,T]$. We also fix $r>0$.

We define $\tilde{X}^\epsilon$ and $\tilde{Y}^{\epsilon}$ via a mild formulation: for any $0\leq t\leq T$
\begin{equation}\label{processhasmild}
\begin{gathered}
\tilde{Y}^\epsilon(t)=e^{\frac{t}{\epsilon}B}y+\frac{1}{\epsilon}\int_{0}^{t}e^{\frac{(t-s)}{\epsilon}B}G(X^{\epsilon}(\lfloor \frac{s}{\delta}\rfloor\delta),\tilde{Y}^{\epsilon}(s)))ds+\frac{1}{\sqrt{\epsilon}}\int_{0}^{t}e^{\frac{(t-s)}{\epsilon}B}dW(s),\\
\tilde{X}^\epsilon(t)=e^{tA}x+\int_{0}^{t}e^{(t-s)A}F(X^\epsilon(\lfloor\frac{s}{\delta}\rfloor\delta),\tilde{Y}^\epsilon(s))ds,
\end{gathered}
\end{equation}
where $\lfloor . \rfloor$ denotes the integer part function.

$\tilde{X}^\epsilon$ and $\tilde{Y}^\epsilon$ are continuous processes and they satisfy $\tilde{X}^\epsilon(0)=x=X^\epsilon(0)$ and $\tilde{Y}^\epsilon(0)=y=Y^\epsilon(0)$.
Moreover, on any subinterval $[k\delta,(k+1)\delta]$, with $0\leq k\leq N:=\lfloor \frac{T}{\delta}\rfloor$, we have
\begin{equation}\label{processhas}
\begin{gathered}
d\tilde{X}^{\epsilon}(t)=(A\tilde{X}^{\epsilon}(t)+F(X^{\epsilon}(k\delta),\tilde{Y}^{\epsilon}(t)))dt,\\
d\tilde{Y}^{\epsilon}(t)=\frac{1}{\epsilon}(B\tilde{Y}^{\epsilon}(t)+G(X^{\epsilon}(k\delta),\tilde{Y}^{\epsilon}(t)))dt+\frac{1}{\sqrt{\epsilon}}dW(t).
\end{gathered}
\end{equation}
We remark that on such subintervals the fast component $\tilde{Y}^\epsilon$ does not depend on the slow component $\tilde{X}^\epsilon$, but only on the value of $X^\epsilon$ at the first point of the interval.

Since $F$ and $G$ are supposed to be bounded, we easily see that for any $t\geq 0$
$$\E|\tilde{X}^\epsilon(t)|_{H}^{2}\leq C(1+|x|_{H}^{2})\quad \text{and} \quad \E|\tilde{Y}^\epsilon(t)|_{H}^{2}\leq C(1+|y|_{H}^{2}).$$

We show the following Lemmas:
\begin{lemme}\label{lemstrong1}
There exists $C>0$ such that for any $0\leq t\leq T$ and any $\epsilon>0$
$$\E|X^\epsilon(t)-\tilde{X}^\epsilon(t)|\leq C(\delta^{1-r}+\frac{\epsilon}{\delta}).$$
\end{lemme}

\begin{lemme}\label{lemstrong2}
There exists $C>0$ such that for any $0\leq t\leq T$ and any $\epsilon>0$
$$\E|\tilde{X}^\epsilon(t)-\overline{X}(t)|\leq
C\left(\epsilon(1+\delta^{-r})(1+\frac{1}{1-e^{-c\frac{\delta}{\epsilon}}})\right)^{1/2}+C(\delta^{1-r}+\frac{\epsilon}{\delta}).$$
\end{lemme}

The constant $C$ above depends on $r$, $T$, $x$, $y$, but not on $t$, $\epsilon$ or $\delta(\epsilon)$.

Lemma \ref{lemstrong1} explains why we can replace $\tilde{X}^\epsilon$ by $X^\epsilon$, and Lemma \ref{lemstrong2} gives an estimate of the distance between $\tilde{X}^\epsilon(t)$ and $\overline{X}(t)$. With the choice of $\delta(\epsilon)=\sqrt{\epsilon}$ given by \eqref{choiceepsilon}, the proof of Theorem \ref{strong} is straightforward.

\underline{Proof of Lemma \ref{lemstrong1}}
\begin{itemize}

\item Estimate of $Y^{\epsilon}-\tilde{Y}^{\epsilon}$.

We define $\rho^{\epsilon}(t)=Y^{\epsilon}(t)-\tilde{Y}^{\epsilon}(t)$ for any $0\leq t\leq T$.

We fix $k\geq 0$; then for any $t\in [k\delta,(k+1)\delta]$ we have
$$
d\rho^{\epsilon}(t)=\frac{1}{\epsilon}B\rho^\epsilon(t)dt
+\frac{1}{\epsilon}(G(X^{\epsilon}(t),Y^{\epsilon}(t))-G(X^{\epsilon}(k\delta),\tilde{Y}^{\epsilon}(t)))dt.
$$

Using a mild formulation and Gronwall Lemma, for any $t\in [k\delta,(k+1)\delta]$ we have
$$\E|\rho^{\epsilon}(t)|\leq e^{-\frac{\mu-L_g}{\epsilon}(t-k\delta)}\E|\rho^\epsilon(k\delta)|
+\frac{C}{\epsilon}\int_{k\delta}^{t}e^{-\frac{\mu-L_g}{\epsilon}(t-s)}\E|X^{\epsilon}(s)-X^{\epsilon}(k\delta)|ds.$$

Since $\E|Y^\epsilon(t)|\leq C(1+|y|)$ and $\E|\tilde{Y}^\epsilon(t)|\leq C(1+|y|)$, we have the same bound on $\rho^\epsilon$.

We can integrate the previous inequality over the interval $t\in[k\delta,(k+1)\delta]$ and get
\begin{align*}
\int_{k\delta}^{(k+1)\delta}\E|\rho^{\epsilon}(t)|dt&\leq C\int_{k\delta}^{(k+1)\delta}e^{-\frac{\mu-L_g}{\epsilon}(t-k\delta)}dt+\frac{C}{\epsilon}\int_{k\delta}^{(k+1)\delta}\int_{k\delta}^{t}e^{-\frac{\mu-L_g}{\epsilon}(t-s)}\E|X^{\epsilon}(s)-X^{\epsilon}(k\delta)|dsdt\\
&\leq C\frac{\epsilon}{\mu-L_g}+C\int_{k\delta}^{(k+1)\delta}\E|X^{\epsilon}(s)-X^{\epsilon}(k\delta)|\int_{s}^{(k+1)\delta}\frac{1}{\epsilon}e^{-\frac{\mu-L_g}{\epsilon}(t-s)}dtds\\
&\leq C\epsilon+C\int_{k\delta}^{(k+1)\delta}\E|X^{\epsilon}(s)-X^{\epsilon}(k\delta)|ds.
\end{align*}

We recall that strict dissipativity $\mu-L_g>0$ holds, thanks to Assumption \ref{strictdiss}.

It remains to take the sum over $k\in\left\{0,\ldots,\lfloor \frac{t}{\delta} \rfloor\right\}$, where $t\leq T$; using Proposition \ref{reg} of the appendix, we then obtain
$$\int_{0}^{t}\E|\rho^{\epsilon}(s)|ds\leq C(r,T)(\frac{\epsilon}{\delta}+\delta^{1-r}).$$

\item Estimate of $X^{\epsilon}-\tilde{X}^{\epsilon}$.

We have for any $0\leq t\leq T$
\begin{align*}
|X^{\epsilon}(t)-\tilde{X}^{\epsilon}(t)|&=|\int_{0}^{t}e^{(t-s)A}(F(X^\epsilon(s),Y^\epsilon(s))-F(X^\epsilon(\lfloor\frac{s}{\delta}\rfloor\delta),\tilde{Y}(s)))ds|\\
&\leq \int_{0}^{t}|F(X^\epsilon(s),Y^\epsilon(s))-F(X^\epsilon(\lfloor\frac{s}{\delta}\rfloor\delta),\tilde{Y}(s))|^2ds\\
&\leq C\int_{0}^{T}(|X^\epsilon(s)-X^\epsilon(\lfloor\frac{s}{\delta}\rfloor\delta)|+|Y^\epsilon(s)-\tilde{Y}^\epsilon(s)|)ds.
\end{align*}

Using the previous estimate and the regularity result from Proposition \ref{reg}, we obtain for any $0\leq t\leq T$
\begin{equation}\label{xxtilde}
\E|X^{\epsilon}(t)-\tilde{X}^{\epsilon}(t)|\leq C(\delta^{1-r}+\frac{\epsilon}{\delta}).
\end{equation}

\end{itemize}
\begin{flushright}
$\Box$
\end{flushright}

\underline{Proof of Lemma \ref{lemstrong2}}
We introduce the following decomposition, for any $0\leq t\leq T$:
\begin{align*}
\tilde{X}^{\epsilon}(t)-\overline{X}(t)&=\int_{0}^{t}e^{(t-s)A}(F(X^{\epsilon}(\lfloor
\frac{s}{\delta}\rfloor
\delta),\tilde{Y}^{\epsilon}(s))-\overline{F}(\overline{X}(s)))ds\\
&=\int_{0}^{t}e^{(t-s)A}(F(X^{\epsilon}(\lfloor \frac{s}{\delta}\rfloor
\delta),\tilde{Y}^{\epsilon}(s))-\overline{F}(X^{\epsilon}(\lfloor
\frac{s}{\delta}\rfloor \delta)))ds\\
&+\int_{0}^{t}e^{(t-s)A}(\overline{F}(X^{\epsilon}(\lfloor \frac{s}{\delta} \rfloor
\delta))-\overline{F}(X^{\epsilon}(s)))ds\\
&+\int_{0}^{t}e^{(t-s)A}(\overline{F}(X^{\epsilon}(s))-\overline{F}(\tilde{X}^{\epsilon}(s)))ds\\
&+\int_{0}^{t}e^{(t-s)A}(\overline{F}(\tilde{X}^{\epsilon}(s))-\overline{F}(\overline{X}(s)))ds\\
&=I_1(t)+I_2(t)+I_3(t)+I_4(t).
\end{align*}

According to Proposition \ref{FbarLip}, $\overline{F}$ is Lipschitz continuous; thanks to Proposition \ref{reg} and Lemma \ref{lemstrong1}, we then show that for any $0\leq t\leq T$
\begin{gather*}
\E|I_2(t)|\leq C\int_{0}^{T}\E|X^{\epsilon}(\lfloor \frac{s}{\delta} \rfloor
\delta)-X^{\epsilon}(s)|ds\leq C\delta^{1-r}\\
\E|I_3(t)|\leq C\int_{0}^{T}|X^{\epsilon}(s)-\tilde{X}^{\epsilon}(s)|ds\leq C(\delta^{1-r}+\frac{\epsilon}{\delta})\\
\E|I_4(t)|\leq CT\int_{0}^{t}\E|\tilde{X}^{\epsilon}(r)-\overline{X}(r)|ds.
\end{gather*}
The $I_4$ term is treated via the Gronwall Lemma.

It remains to focus on the $I_1$ term. It is fundamental to look at $\E|I_1(t)|^2$ and not only $\E|I_1(t)|$ in order to obtain the best estimate leading to order $1/2$. For that, we use the subdivision of $[0,T]$ and we expand the scalar product in $H$: we have for any $0\leq t\leq T$
\begin{align*}
|I_1(t)|^2&=|\int_{0}^{t}e^{(t-s)A}(F(X^{\epsilon}(\lfloor\frac{s}{\delta}\rfloor\delta),\tilde{Y}^{\epsilon}(s))-\overline{F}(X^{\epsilon}(\lfloor\frac{s}{\delta}\rfloor\delta)))ds|^2\\
&=|\sum_{k=0}^{\lfloor\frac{t}{\delta}\rfloor}\int_{k\delta}^{(k+1)\delta\wedge t}e^{(t-s)A}(F(X^{\epsilon}(k\delta),\tilde{Y}^{\epsilon}(s))-\overline{F}(X^{\epsilon}(k\delta)))ds|^2\\
&=A_1(t)+A_2(t),
\end{align*}
where
\begin{equation}\label{defA1}
A_1(t):=\sum_{k=0}^{\lfloor\frac{t}{\delta}\rfloor}|\int_{k\delta}^{(k+1)\delta\wedge t}e^{(t-s)A}(F(X^{\epsilon}(k\delta),\tilde{Y}^{\epsilon}(s))-\overline{F}(X^{\epsilon}(k\delta)))ds|^2
\end{equation}
and
\begin{multline}\label{defA2}
A_2(t):=2\sum_{0\leq i<j\leq \lfloor\frac{t}{\delta}\rfloor}\langle\int_{i\delta}^{(i+1)\delta\wedge t}e^{(t-s)A}(F(X^{\epsilon}(i\delta),\tilde{Y}^{\epsilon}(s))-\overline{F}(X^{\epsilon}(i\delta)))ds,\\
\int_{j\delta}^{(j+1)\delta\wedge t}e^{(t-s)A}(F(X^{\epsilon}(j\delta),\tilde{Y}^{\epsilon}(s))-\overline{F}(X^{\epsilon}(j\delta)))ds\rangle.
\end{multline}

We claim that $\E A_1(t)\leq C\epsilon$ and $\E A_2(t)\leq C\epsilon(1+\delta^{-r})(1+\frac{1}{1-e^{-c\frac{\delta}{\epsilon}}})$, where $C=C(T,r,\theta,x,y)$. Using Gronwall Lemma, we get the result.

\begin{itemize}
\item We first prove the estimate on $\E A_1(t)$. We use conditional expectation with respect to $\mathcal{F}_s$. Then for any $0\leq k\leq\lfloor\frac{t}{\delta}\rfloor$, using some symmetry for variables $s$ and $\sigma$,
\begin{align*}
\E&|\int_{k\delta}^{(k+1)\delta\wedge t}e^{(t-s)A}(F(X^{\epsilon}(k\delta),\tilde{Y}^{\epsilon}(s))-\overline{F}(X^{\epsilon}(k\delta)))ds|^2\\
&=2\E\int_{k\delta}^{(k+1)\delta\wedge t}ds\int_{s}^{(k+1)\delta\wedge t}d\sigma\langle e^{(t-s)A}(F(X^{\epsilon}(k\delta),\tilde{Y}^{\epsilon}(s))-\overline{F}(X^{\epsilon}(k\delta))),\\
&\hspace{150 pt}e^{(t-\sigma)A}(F(X^{\epsilon}(k\delta),\tilde{Y}^{\epsilon}(\sigma))-\overline{F}(X^{\epsilon}(k\delta)))\rangle \\
&=2\E\int_{k\delta}^{(k+1)\delta\wedge t}ds\int_{s}^{(k+1)\delta\wedge t}d\sigma\langle e^{(t-s)A}(F(X^{\epsilon}(k\delta),\tilde{Y}^{\epsilon}(s))-\overline{F}(X^{\epsilon}(k\delta))),\\
&\hspace{150 pt}e^{(t-\sigma)A}\E[F(X^{\epsilon}(k\delta),\tilde{Y}^{\epsilon}(\sigma))-\overline{F}(X^{\epsilon}(k\delta))|\mathcal{F}_s]\rangle\\
&\leq 2\int_{k\delta}^{(k+1)\delta}\int_{s}^{(k+1)\delta}\E\left(|F(X^{\epsilon}(k\delta),\tilde{Y}^{\epsilon}(s))-\overline{F}(X^{\epsilon}(k\delta))||\E[F(X^{\epsilon}(k\delta),\tilde{Y}^{\epsilon}(\sigma))|\mathcal{F}_s]-\overline{F}(X^{\epsilon}(k\delta))|\right).
\end{align*}

We now define the auxiliary function $\tilde{F}$. Propositions \ref{tildeFinf} and \ref{tildeFLip} - see the Appendix - give important properties of this function: we have some exponential control with respect to time $t$ of uniform and Lipschitz bounds with respect to $x$.
\begin{defi}\label{defiauxfunctions}
For any $(x,y)\in H^2$ and $t\geq 0$
\begin{equation}\label{defFtilde}
\tilde{F}(x,y,t)=\E F(x,Y_{x}(t,y))-\overline{F}(x).
\end{equation}
\end{defi}

Since $F$ is bounded, we can control the first factor in the integral; for the second factor, we use the definition of $\tilde{F}$, the Markov property and Proposition \ref{tildeFinf} to see that
\begin{align*}
\E|&\int_{k\delta}^{(k+1)\delta\wedge t}e^{(t-s)A}(F(X^{\epsilon}(k\delta),\tilde{Y}^{\epsilon}(s))-\overline{F}(X^{\epsilon}(k\delta)))ds|^2\\
&\leq C\int_{k\delta}^{(k+1)\delta}ds\int_{s}^{(k+1)\delta}d\sigma\E|\tilde{F}(X^\epsilon(k\delta),\tilde{Y}^\epsilon(s),\frac{\sigma-s}{\epsilon})|\\
&\leq C\int_{k\delta}^{(k+1)\delta}\int_{s}^{(k+1)\delta}e^{-c\frac{\sigma-s}{\epsilon}}d\sigma ds\leq C\delta\epsilon.
\end{align*}
Therefore we get $\E A_1(t)\leq C\epsilon$.

\item Estimate of $\E A_2(t)$.

We have to introduce the following auxiliary processes, which generalize $\tilde{Y}^\epsilon$.
$(Z_i^{\epsilon}(s))_{s\geq i\delta}$, where $i\in\{0,\ldots,N\}$ is defined by:
\begin{equation}\label{procZ}
\begin{gathered}
dZ_{i}^{\epsilon}(s)=\frac{1}{\epsilon}(BZ_{i}^{\epsilon}(s)+G(X^\epsilon(i\delta),Z_{i}^{\epsilon}(s)))ds+\frac{1}{\sqrt{\epsilon}}dW(s)\\
Z_{i}^{\epsilon}(i\delta)=\tilde{Y}^\epsilon(i\delta).
\end{gathered}
\end{equation}
It is then clear that for $i\delta\leq s\leq (i+1)\delta$ we have $Z_{i}^{\epsilon}(s)=\tilde{Y}^\epsilon(s)$, and that $Z_{k}^{\epsilon}((k+1)\delta)=\tilde{Y}^\epsilon((k+1)\delta)=Z_{k+1}^{\epsilon}((k+1)\delta)$. Moreover the processes $(Z_{i}^{\epsilon})$ are uniformly bounded with respect to $i$ and $\epsilon$.

It is then possible to rewrite the integrands appearing in the expression of $A_2$: when $i\delta\leq s\leq (i+1)\delta\leq j\delta\leq \tau\leq (j+1)\delta$,
\begin{align*}
&|\E\left\langle e^{(t-s)A}(F(X^\epsilon(i\delta),\tilde{Y}^\epsilon(s))-\overline{F}(X^\epsilon(i\delta))),e^{(t-\tau)A}(F(X^\epsilon(j\delta),\tilde{Y}^\epsilon(\tau))-\overline{F}(X^\epsilon(j\delta)))\right\rangle|\\
&=|\E\left\langle e^{(t-s)A}(F(X^\epsilon(i\delta),\tilde{Y}^\epsilon(s))-\overline{F}(X^\epsilon(i\delta))),e^{(t-\tau)A}\E[F(X^\epsilon(j\delta),\tilde{Y}^\epsilon(\tau))-\overline{F}(X^\epsilon(j\delta))|\mathcal{F}_{(i+1)\delta}]\right\rangle|\\
&\leq C\E|\E[F(X^\epsilon(j\delta),\tilde{Y}^\epsilon(\tau))-\overline{F}(X^\epsilon(j\delta))|\mathcal{F}_{(i+1)\delta}]|,
\end{align*}
since $F$ is bounded.
We have $\tilde{Y}^\epsilon(\tau)=Z_{j}^{\epsilon}(\tau)$; however since $i<j$ we can use conditional expectation with respect to $\mathcal{F}_{(i+1)\delta}$ instead of $\mathcal{F}_{j\delta}$  in order to get a better estimate. We therefore propose the following decomposition
\begin{align*}
&\E|\E[F(X^\epsilon(j\delta),\tilde{Y}^\epsilon(\tau))-\overline{F}(X^\epsilon(j\delta))|\mathcal{F}_{(i+1)\delta}]|\\
&\leq C\E|\E[F(X^\epsilon((i+1)\delta),Z_{i+1}^{\epsilon}(\tau))-\overline{F}(X^\epsilon((i+1)\delta))|\mathcal{F}_{(i+1)\delta}]|\\
&+C\E|\E\left[(F(X^\epsilon(j\delta),Z_{j}^{\epsilon}(\tau))-\overline{F}(X^\epsilon(j\delta)))-(F(X^\epsilon((i+1)\delta),Z_{i+1}^{\epsilon}(\tau))-\overline{F}(X^\epsilon((i+1)\delta)))|\mathcal{F}_{(i+1)\delta}\right]|\\
&=:B_1+B_2.
\end{align*}

\begin{enumerate}
\item First, using Markov property we have for any $j\delta\leq \tau\leq (j+1)\delta$
\begin{align*}
B_1&=C\E|\E\left[F(X^\epsilon((i+1)\delta),Z_{i+1}^{\epsilon}(\tau))-\overline{F}(X^\epsilon((i+1)\delta))|\mathcal{F}_{(i+1)\delta}\right]|\\
&=C\E|\E\left[\tilde{F}(X^\epsilon((i+1)\delta),Z_{(i+1)}^{\epsilon}((i+1)\delta),\frac{\tau-(i+1)\delta}{\epsilon})|\mathcal{F}_{(i+1)\delta}\right]|\\
&\leq Ce^{-c\frac{\tau-(i+1)\delta}{\epsilon}},
\end{align*}
thanks to Proposition \ref{tildeFinf}.
\item $B_2$ can be rewritten using a telescoping sum, and some conditional expectation
\begin{align*}
B_2&=\E|\E\Big(\sum_{k=i+1}^{j-1}\E\Big[[F(X^\epsilon((k+1)\delta),Z_{k+1}^{\epsilon}(\tau))-\overline{F}(X^\epsilon((k+1)\delta))]\\
&\hspace{200 pt}-[F(X^\epsilon(k\delta),Z_{k}^{\epsilon}(\tau))-\overline{F}(X^\epsilon(k\delta))]|\mathcal{F}_{k\delta}\Big]|\mathcal{F}_{(i+1)\delta}\Big)|\\
&\leq \sum_{k=i+1}^{j-1}\E|\E\Big[[F(X^\epsilon((k+1)\delta),Z_{k+1}^{\epsilon}(\tau))-\overline{F}(X^\epsilon((k+1)\delta))]\\
&\hspace{200 pt}-[F(X^\epsilon(k\delta),Z_{k}^{\epsilon}(\tau))-\overline{F}(X^\epsilon(k\delta))]|\mathcal{F}_{k\delta}\Big]|.
\end{align*}
Thanks to Markov property, we obtain
\begin{align*}
\E\Big[F(X^\epsilon((k+1)\delta),Z_{k+1}^{\epsilon}(\tau))&-\overline{F}(X^\epsilon((k+1)\delta))|\mathcal{F}_{k\delta}\Big]\\
&=\E\Big[\tilde{F}(X^\epsilon((k+1)\delta),\tilde{Y}^\epsilon((k+1)\delta),\tau-(k+1)\delta,\epsilon)|\mathcal{F}_{k\delta}\Big]
\end{align*}
and
$$
\E\Big[F(X^\epsilon(k\delta),Z_{k}^{\epsilon}(\tau))-\overline{F}(X^\epsilon(k\delta)))|\mathcal{F}_{k\delta}\Big]=\E\Big[\tilde{F}(X^\epsilon(k\delta),\tilde{Y}^\epsilon((k+1)\delta),\tau-(k+1)\delta,\epsilon)|\mathcal{F}_{k\delta}\Big].
$$

Using the exponential decrease in time of the Lipschitz constant of $\tilde{F}$ with respect to $x$ given by Proposition \ref{tildeFLip} in the appendix, we have
\begin{align*}
B_2&\leq \sum_{k=i+1}^{j-1}\E|\E\Big[\tilde{F}(X^\epsilon((k+1)\delta),\tilde{Y}^\epsilon((k+1)\delta),\tau-(k+1)\delta,\epsilon)|\mathcal{F}_{k\delta}\Big]\\
&\hspace{100 pt}-\E\Big[\tilde{F}(X^\epsilon(k\delta),\tilde{Y}^\epsilon((k+1)\delta),\tau-(k+1)\delta,\epsilon)|\mathcal{F}_{k\delta}\Big]|\\
&\leq C\sum_{k=i+1}^{j-1}e^{-c\frac{\tau-(k+1)\delta}{\epsilon}}\E|X^\epsilon((k+1)\delta)-X^\epsilon(k\delta)|(1+\frac{\epsilon^\eta}{(\tau-(k+1)\delta)^\eta})\\
&\leq C(1+\frac{\epsilon^\eta}{(\tau-j\delta)^\eta})\sum_{k=i+1}^{j-1}e^{-c\frac{\tau-(k+1)\delta}{\epsilon}}\frac{\delta^{1-r}}{(k\delta)^{1-r}}\\
&\leq C(1+\frac{\epsilon^\eta}{(\tau-j\delta)^\eta})\frac{\delta^{1-r}}{((i+1)\delta)^{1-r}}\frac{e^{-c\frac{\tau-j\delta}{\epsilon}}}{1-e^{-c\frac{\delta}{\epsilon}}},
\end{align*}
by using the regularity proved in Proposition \ref{reg}.

\item We are now able to conclude, by using the estimates on $B_1$ and $B_2$: we have for any $0\leq t\leq T$
\begin{align*}
\E A_2(t)&\leq C\sum_{0\leq i<j\leq \lfloor\frac{t}{\delta}\rfloor}\int_{i\delta}^{(i+1)\delta\wedge t}ds\int_{j\delta}^{(j+1)\delta\wedge t}d\tau e^{-c\frac{\tau-(i+1)\delta}{\epsilon}}\\
&+C\sum_{0\leq i<j\leq \lfloor\frac{t}{\delta}\rfloor}\int_{i\delta}^{(i+1)\delta\wedge t}ds\int_{j\delta}^{(j+1)\delta\wedge t}d\tau \frac{\delta^{1-r}}{((i+1)\delta)^{1-r}}\frac{e^{-c\frac{\tau-j\delta}{\epsilon}}}{1-e^{-c\frac{\delta}{\epsilon}}}(1+\frac{\epsilon^\eta}{(\tau-j\delta)^\eta})\\
&\leq C\sum_{0\leq i<j\leq \lfloor\frac{T}{\delta}\rfloor}\epsilon\delta e^{-c\frac{(j-i-1)\delta}{\epsilon}}+C\delta^{1-r}\epsilon\frac{1}{1-e^{-c\frac{\delta}{\epsilon}}}\sum_{0\leq i<j\leq \lfloor\frac{T}{\delta}\rfloor}\frac{\delta^{1-r}}{((i+1)\delta)^{1-r}}\\
&\leq C(T)\epsilon(1+\delta^{-r})(1+\frac{1}{1-e^{-c\frac{\delta}{\epsilon}}}).
\end{align*}

\end{enumerate}

\end{itemize}

\begin{flushright}
$\Box$
\end{flushright}

\section{Proof of the weak-order result}\label{weakproofsect}

The proof of Theorem \ref{weak} relies on an expansion of the solution of the Kolmogorov equation associated with the stochastic system \eqref{eqgen} with respect to the small parameter $\epsilon$; the zero-order term corresponds to the averaged equation, and we control the first-order term to get the result.

When working with SDEs, this strategy can be entirely followed; nevertheless in the case of SPDEs, the Kolmogorov equations involve the unbounded operators $A$ and $B$, and quantities like $\text{Tr}(D_{yy}^{2}u)$, and this leads to technical problems - see \cite{CeBook} or \cite{Ce-F}.

The idea of the proof is to reduce the infinite dimensional problem to a finite dimensional one by Galerkin approximation; we apply the method in this finite dimensional setting, and we prove bounds that are uniform with respect to the dimension. We also show that taking the limit when dimension goes to infinity is meaningful and gives the desired result.

The key element in the construction of the expansion mentioned above is given in Lemma \ref{lemPoiss} below: a Poisson equation can be solved under ergodicity conditions on the fast equation. We notice that Assumption \ref{weakdiss} is sufficient, since we can analyse the problem through the asymptotic properties of the transition semi-group of the fast equation, instead of trajectories.


First, we explain how we reduce the problem to a finite dimensional one - see Section \ref{sectionreduc}; then we explain the method in this setting and show which expressions must be controlled - see Section \ref{sectionexpansion}; finally we prove the estimates.

\subsection{Reduction to a finite dimensional problem}\label{sectionreduc}

We use Galerkin approximations based on the orthonormal basis $(e_k)$ and $(f_k)$ of $H$, given by Assumption \ref{hypAB}. We define the subspaces
$$
H_{N}^{(1)}=\text{span}\left\{ e_k; 0\leq k\leq N-1\right\} \quad \text{and} \quad H_{N}^{(2)}=\text{span}\left\{ f_k; 0\leq k\leq N-1\right\}.
$$
We denote by $P_{N}^{(1)}\in\mathcal{L}(H)$ - resp. $P_{N}^{(2)}$ - the orthogonal projection of $H$ onto $H_{N}^{(1)}$ - resp. $H_{N}^{(2)}$.

Then we define for $x\in H_{N}^{(1)}$ and $y\in H_{N}^{(2)}$
\begin{gather*}
F_{N}(x,y)=P_{N}^{(1)}(F(x,y))\\
G_{N}(x,y)=P_{N}^{(2)}(G(x,y))\\
U_{N}(x,y)=U(x,y).
\end{gather*}
On $H_{N}^{(1)}\times H_{N}^{(2)}$ coefficients $F_N$ and $G_N$ are of class $\mathcal{C}^2$. The function $U_N$ is of class $\mathcal{C}^3$ with respect to $y$ and of class $\mathcal{C}^2$ with respect to $x$, and we have $D_yU_N(x,y)=G_N(x,y)$ for any $(x,y)\in H_{N}^{(1)}\times H_{N}^{(2)}$.

Moreover we have bounds on $F_N$, $G_N$, $U_N$ and on their derivatives which are uniform with respect to $N$ and satisfy bound like in Assumption \ref{hypF}. In particular we still have the weak dissipativity condition with $G$ replaced by $G_N$.

We can define the following approximation of system \eqref{eqgen}
\begin{equation}\label{eqproj}
\begin{gathered}
dX_{N}^{\epsilon}(t)=(AX_{N}^{\epsilon}(t)+F_N(X_{N}^{\epsilon}(t),Y_{N}^{\epsilon}(t)))dt\\
dY_{N}^{\epsilon}(t)=\frac{1}{\epsilon}(BY_{N}^{\epsilon}(t)+G_{N}(X_{N}^{\epsilon}(t),Y_{N}^{\epsilon}(t)))dt+\frac{1}{\sqrt{\epsilon}}dW_{N}(t),
\end{gathered}
\end{equation}
with initial conditions $X_{N}^{\epsilon}(0)=P_{N}^{(1)}x\in H_{N}^{(1)}$, $Y_{N}^{\epsilon}(0)=P_{N}^{(2)}y\in H_{N}^{(2)}$.

Since $H_{N}^{(1)}$ (resp. $H_{N}^{(2)}$) is a stable subspace of $A$ (resp. $B$), this system is well-posed in $H_{N}^{(1)}\times H_{N}^{(2)}$.

We have by definition $W_{N}(t)=P_{N}^{(2)}W(t)$; on $H_{N}^{(2)}$ it is a $N$-dimensional Brownian motion.

Below we explain that system \eqref{eqproj} defines a good approximation of the initial problem \eqref{eqgen}. Moreover, we can check that the structure of the problem remains the same, with bounds independent of the dimension.

First we describe the ergodic properties, and in particular the relations between the invariant measures associated with the associated fast equations with frozen slow component. In Section \ref{known}, we have defined $\nu=\mathcal{N}(0,\frac{(-B)^{-1}}{2})$ and $\mu^x$ - see \eqref{definvmeas}; we can do the same for the finite dimensional fast equation with frozen slow component $x\in H$ - and not only for $x\in H_{N}^{(1)}$:
\begin{equation}\label{eqprojfig}
\begin{gathered}
dY_{x,N}(t,y)=(BY_{x,N}(t,y)+P_{N}^{(2)}G(x,P_{N}^{(2)}Y_{x,N}(t,y)))dt+dW_N(t)\\
Y_{x,N}(0,y)=P_{N}^{(2)}y.
\end{gathered}
\end{equation}

Let $\nu_N$ be the unique centered Gaussian probability measure on $H_{N}^{(2)}$ having for covariance operator the induced matrix from $\frac{(B)^{-1}}{2}$ on the subspace $H_{N}^{(2)}$. We can then build $\mu_{N}^{x}$ the unique - since we have strict dissipativity - invariant probability measure associated to \eqref{eqprojfig}: we naturally extend the definition of $U_N$ to $H\times H_{N}^{(2)}$, by $U_N(x,y)=U(x,y)$, and we define
\begin{equation}\label{mesfinie}
\mu_{N}^{x}(dy)=\frac{1}{Z_N(x)}e^{2U_N(x,y)}\nu_N(dy),
\end{equation}
where $Z_N(x)\in]0,+\infty[$ is a normalization constant.

As $\nu_N$ is the image measure of $\nu$ by $P_{N}^{(2)}$, we can use the following change of variables formula for suitable test functions $\Phi$:
\begin{align*}
\int_{H_{N}^{(2)}}\Phi(y)\mu_{N}^{x}(dy)&=\frac{1}{Z_N(x)}\int_{H_{N}^{(2)}}\Phi(y)e^{2U(x,y)}\nu_N(dy)\\
&=\frac{1}{Z_N(x)}\int_{H}\Phi(P_{N}^{(2)}y)e^{2U(x,P_{N}^{(2)}y)}\nu(dy).
\end{align*}
This formula leads to convergence properties: first when $N\rightarrow +\infty$ we have
$$Z_N(x)\rightarrow Z(x) \text{ for any }x\in H.$$
Moreover if $\Phi:H\rightarrow H$ is a continuous function such that for any $y\in H$, $|\Phi(y)|\leq c(1+|y|)$, then
$$\int_{H_{N}^{(2)}}\Phi(y)\mu_{N}^{x}(dy)\rightarrow \int_{H}\Phi(y)\mu^x(dy).$$

For any $x\in H$ we define the averaged coefficient associated with the finite dimensional problem \eqref{eqprojfig}
\begin{equation}\label{coeffmoyfinie}
\overline{F_{N}}(x)=\int_{H_{N}^{(2)}}P_{N}^{(1)}F(x,y)\mu_{N}^{x}(dy),
\end{equation}
and the new averaged equation
\begin{equation}\label{eqprojmoy}
d\overline{X_{N}}(t)=(A\overline{X_{N}}(t)+\overline{F_N}(\overline{X_{N}}(t)))dt,
\end{equation}
with initial condition $\overline{X_{N}}(0)=P_{N}^{(1)}x\in H_{N}^{(1)}$.

We notice that $\overline{F_N}$ is bounded, and of class $\mathcal{C}^2$ with bounded derivatives; moreover it satifies the properties described in Proposition \ref{FbarLip}, with constants independent of $N$.

\begin{rem}
Making the Galerkin projection and then averaging the coefficient with $\mu_{N}^{x}$ is not the same as averaging the coefficient with $\mu^{x}$ and then making the Galerkin projection. As a consequence $\overline{X_{N}}$ is not naturally defined by a Galerkin approximation from $\overline{X}$.
\end{rem}

The following Lemma gives the convergence of the finite dimensional approximations to the initial problem:
\begin{lemme}\label{finieversinfinie}
\begin{enumerate}
 \item For any fixed $\epsilon>0$, $t\geq 0$, and any $x\in H$, $y\in H$, we have when $N\rightarrow +\infty$
$$\E|X^{\epsilon}(t)-X_{N}^{\epsilon}(t)|^2+\E|Y^{\epsilon}(t)-Y_{N}^{\epsilon}(t)|^2\rightarrow 0.$$
 \item For any $t\geq 0$, $x\in H$, we have when $N\rightarrow +\infty$
$$|\overline{X}(t)-\overline{X_{N}}(t)|\rightarrow0.$$
\end{enumerate}
\end{lemme}

It remains to define an approximated test function: for any $x\in H_{N}^{(1)}$, we define $\phi_{N}(x)=\phi(x)$.

We can now show how the initial problem can be reduced to a finite dimensional one - provided we are able to give estimates uniform to the dimension: for any initial conditions $x,y\in H$ and $0\leq t\leq T$,
\begin{align*}
\E[\phi(X^{\epsilon}(t))]-\E[\phi(\overline{X}(t))]&=\E[\phi(X^{\epsilon}(t))]-\E[\phi(X_{N}^{\epsilon}(t))]\\
&+\E[\phi_{N}(X_{N}^{\epsilon}(t)]-\E[\phi_{N}(\overline{X_{N}}(t))]\\
&+\E[\phi(\overline{X_{N}}(t))]-\E[\phi(\overline{X}(t))].
\end{align*}

The first and the last terms converge to $0$ when $N\rightarrow +\infty$, according to the above Lemma \ref{finieversinfinie}; we later control the central term with an expression independent of $N$. Taking the limit as dimension goes to infinity then gives the result. So we have to control
\begin{equation}\label{weakN}
\E[\phi_{N}(X_{N}^{\epsilon}(t)]-\E[\phi_{N}(\overline{X_{N}}(t))].
\end{equation}
From now on, we only work with the approximations to obtain estimates, but in order to simplify the notations, we forget the index $N$ - since bounds are uniform with respect to $N$. The spaces $H_{N}^{(i)}$ are denoted by $H^{(i)}$ in the next sections.

\subsection{The asymptotic expansion}\label{sectionexpansion}

We define the following differential operators: for functions of class $\mathcal{C}^2$ $\psi:H^{(1)}\times H^{(2)}\rightarrow \R$, for any $(x,y)\in H^{(1)}\times H^{(2)}$
\begin{gather*}
L_1\psi(x,y)=<By+G(x,y),D_{y}\psi(x,y)>+\frac{1}{2}\text{Tr}(D_{yy}^{2}\psi(x,y))\\
L_2\psi(x,y)=<Ax+F(x,y),D_{x}\psi(x,y)>\\
L^{\epsilon}=\frac{1}{\epsilon}L_1+L_2.
\end{gather*}

We also define for $\psi:H^{(1)}\rightarrow \R$ of class $\mathcal{C}^1$ $\overline{L}\psi(x)=<Ax+\overline{F}(x),D_{x}\psi(x)>$.

We define the following functions $u^\epsilon$ and $\overline{u}$: for $x\in H^{(1)}$, $y\in H^{(2)}$ and $t\geq 0$
\begin{equation}\label{defu}
\begin{gathered}
u^{\epsilon}(t,x,y)=\E[\phi(X^{\epsilon}(t,x,y))]\\
\overline{u}(t,x)=\phi(\overline{X}(t,x)),
\end{gathered}
\end{equation}
where we have mentioned explicitly the dependence on the initial conditions $x,y$ in $X^\epsilon$ and $\overline{X}$.

Since the test function $\phi$ is of class $\mathcal{C}_{b}^{2}$, $u^\epsilon$ and $\overline{u}$ are of class $\mathcal{C}^1$ with respect to $t$ and of class $\mathcal{C}_{b}^{2}$ with respect to $x,y$; we also know that $u^\epsilon$ and $\overline{u}$ are solutions of the following Kolmogorov equations:
\begin{equation}
\begin{gathered}
\frac{\partial u^{\epsilon}}{\partial t}(t,x)=L^{\epsilon}u^{\epsilon}(t,x)\\
u^{\epsilon}(0,x)=\phi(x);
\end{gathered}
\end{equation}

\begin{equation}
\begin{gathered}
\frac{\partial \overline{u}}{\partial t}(t,x,y)=\overline{L}\overline{u}(t,x,y)\\
\overline{u}(0,x,y)=\phi(x).
\end{gathered}
\end{equation}
We remark that the second equation is a linear transport equation with no diffusion term. We then rewrite the expression we want to study (see \eqref{weakN}):
\begin{equation}\label{remkol}
\E[\phi(X^{\epsilon}(T,x,y))]-\E[\phi(\overline{X}(T,x))]=u^{\epsilon}(T,x,y)-\overline{u}(T,x).
\end{equation}

Our strategy is to look for an expansion of $u^{\epsilon}$ with respect to the small parameter $\epsilon$:
\begin{equation}\label{expuepsilon}
u^{\epsilon}=u_{0}+\epsilon u_{1}+v^\epsilon,
\end{equation}
where $v^\epsilon$ is a residual term, while $u_0$ and $u_1$ are smooth and are constructed below.

The identification with respect to the powers of $\epsilon$ gives the following equations:
\begin{equation}\label{cons1}
L_{1}u_0=0,
\end{equation}

\begin{equation}\label{cons2}
\frac{\partial u_0}{\partial t}=L_1u_1+L_2u_0.
\end{equation}


The operator $L_1$ satisfies the following property on the solutions of Poisson equations:
\begin{lemme}\label{lemPoiss}
We fix $x\in H$.
\begin{itemize}
\item If $\Psi$ is a bounded continuous function such that $\int_{H}\Psi(y)\mu^x(dy)=0$, then if $\Phi$ is a function of class $\mathcal{C}^2$ satisfying $L_1\Phi=-\Psi$ then for any $y\in H$ we have
$$\Phi(y)=\int_{H}\Phi\mu^x+\int_{0}^{+\infty}\E[\Psi(Y_x(s,y))]ds.$$
\item Moreover if $\Psi$ is of class $\mathcal{C}_{b}^{2}$, then $\Phi$ defined by
$$\Phi(y)=\int_{0}^{+\infty}\E[\Psi(Y_x(s,y))]ds$$
is of class $\mathcal{C}^2$, satifies $L_1\Phi=-\Psi$, and there exists a constant $C$ - independent on $N$ - such that for any $y\in H$ we have
$$|\Phi(y)|\leq C(1+|y|^2)\|\Psi\|_{\infty}.$$
\end{itemize}
\end{lemme}

\underline{Proof}
The first part of the Lemma is an easy consequence of It\^o formula and of equation \eqref{cvexpy1y2}, after integration with respect to $y_2$ under $\mu^x$. To prove the second part of the lemma, we first see that for any fixed $s\in \R^+$ the function $y\mapsto \E[\Psi(Y_x(s,y))]=:v_x(s,y)$ is of class $\mathcal{C}^2$. To be able to exchange integration in $s$ and derivation with respect to $y$, we need to prove an estimate of the first and the second derivatives which is integrable with respect to $s$.
The derivatives of $Y_x(s,y)$ with respect to the initial condition satisfy the following equations - to simplify notations we do not write dependence in $x$ in those derivatives:
\begin{gather*}
\frac{d\tilde{\eta}^{h,y}(s)}{ds}=B\tilde{\eta}^{h,y}(s)+D_yG(x,Y_x(s,y)).\tilde{\eta}^{h,y}(s,y)\\
\tilde{\eta}^{h,y}(0)=h,
\end{gather*}
and
\begin{gather*}
\frac{d\tilde{\zeta}^{h,k,y}(s)}{ds}=B\tilde{\zeta}^{h,k,y}(s)+D_yG(x,Y_x(s,y)).\tilde{\zeta}^{h,k,y}(s)+D_{yy}^{2}G(x,Y_x(s,y)).(\tilde{\eta}^{h,y}(s),\tilde{\eta}^{k,y}(s))\\
\tilde{\zeta}^{h,k,y}(0)=0.
\end{gather*}

Without any further dissipativity assumption than \eqref{hypdiss}, we only get bounds on finite time intervals like $[0,1]$: there exists $C>0$ such that for any $y,h,k\in H$ and $0\leq s\leq 1$
\begin{gather*}
|D_{y}v_x(s,y).h|\leq C|h|\\
|D_{yy}^2v_x(s,y).(h,k)|\leq C|h||k|.
\end{gather*}

However, using the estimate \eqref{cvexpy1y2} and a Bismut-Elworthy-Li formula, we can indeed prove some exponential convergence with respect to $s$ of the derivatives $D_yv_x(s,y)$ and $D_{yy}^{2}v_x(s,y)$.

Let $\Psi^0$ be a function such that $|\Psi^0(y)|\leq C(\Psi^0)(1+|y|^2)$ for any $y\in H$. If we define $v_{x}^{0}(s,y):=\E\Psi^0(Y_x(s,y))$ for any $y\in H$, $h,k\in H$, we get for the first order derivative
\begin{equation}\label{BEL1}
\begin{aligned}
D_{y}v_{x}^{0}(s,y).h&=\frac{1}{s}\E[\int_{0}^{s}<\tilde{\eta}^{h,y}(\sigma,dW(\sigma)>\Psi^0(Y_x(s,y))]\\
&=\frac{2}{s}\E[\int_{0}^{s/2}<\tilde{\eta}^{h,y}(\sigma),dW(\sigma)>v_{x}^{0}(s/2,Y_x(s/2,y))],
\end{aligned}
\end{equation}
with the observation that $v_{x}^{0}(s,y)=\E v_{x}^{0}(s/2,Y_x(s/2,y))$ thanks to the Markov property; the second order derivative satisfies
\begin{equation}\label{BEL2}
\begin{aligned}
D_{yy}^2v_{x}^{0}(s,y).(h,k)&=\frac{2}{s}\E[\int_{0}^{s/2}\tilde{\zeta}^{h,k,y}(\sigma),dW(\sigma)>v_{x}^{0}(s/2,Y(s/2,y))]\\
&+\frac{2}{s}\E[\int_{0}^{s/2}<\tilde{\eta}^{h,y}(\sigma),dW(\sigma)>D_{y}v_{x}^{0}(s/2,Y(s/2)).\tilde{\eta}^{k,y}(s/2)].
\end{aligned}
\end{equation}
Since $Y_x$ can be controlled in a $L^2$ norm according to Lemma \ref{aprioribounds} in the appendix, we see that there exists $C>0$ such that for any $0<s\leq 1$, $y\in H$, $h,k\in H$
\begin{equation}\label{consBEL}
\begin{gathered}
|D_{y}v_{x}^{0}(s,y).h|\leq \frac{C}{\sqrt{s}}C(\Psi^0)(1+|y|^2)|h|,\\
|D_{yy}^2v_{x}^{0}(s,y).(h,k)||\leq \frac{C}{s}C(\Psi^0)(1+|y|^2)|h||k|.
\end{gathered}
\end{equation}

Now when $s\geq 1$ the Markov property implies that $v_x(s,y)=\E v_x(s-1,Y_x(1,y))$, and choosing $y_1=y$ and by integrating with respect to $\mu^x(dy_2)$ in \eqref{cvexpy1y2} we have
$$|v_x(s-1,y)|\leq Ce^{-c(s-1)}(1+|y|^2).$$
By \eqref{consBEL} at time $1$, we obtain for $s\geq 1$
\begin{gather*}
|D_{y}v_x(s,y).h|\leq Ce^{-c(s-1)}(1+|y|^2)|h|\\
|D_{yy}^2v_x(s,y).(h,k)|\leq Ce^{-c(s-1)}(1+|y|^2)|h||k|.
\end{gather*}
Moreover we have a uniform control when $0\leq s\leq 1$, so that with a change of constants we get the result.

\begin{flushright}
$\Box$
\end{flushright}


As a consequence we see from \eqref{cons1} that $u_0$ is independent of $y$; we then write $u_0(t,x,y)=u_0(t,x)$. We also choose the initial condition $u_0(0,x)=\phi(x)$. The second equation \eqref{cons2} then yields

\begin{align*}
\frac{\partial u_0}{\partial t}(t,x)&=\int_{H^{(2)}}\frac{\partial u_0}{\partial t}(t,x)\mu^{x}(dy)\\
&=\int_{H^{(2)}}L_{1}u_1(t,x,y)\mu^{x}(dy)+\int_{H^{(2)}}L_2u_0(t,x)\mu^{x}(dy)\\
&=<Au_0(t,x)+\int_{H^{(2)}}F(x,y)\mu^{x}(dy),D_{x}u_0(t,x)>\\
&=\overline{L}u_0(t,x).
\end{align*}

\noindent $u_0$ and $\overline{u}$ are solutions of the same evolution equation, with the same initial condition; we can then conclude that $u_0=\overline{u}$.

Then the second equation can be transformed: $\overline{L}u_0=L_1u_1+L_2u_0$. We then obtain an equation on $u_1$:
\begin{equation}\label{defchi}
L_1u_1(t,x,y)=<\overline{F}(x)-F(x,y),D_{x}u_0(t,x)>=:-\chi(t,x,y),
\end{equation}

where $\chi$ is of class $\mathcal{C}_{b}^{2}$ with respect to $y$, and satisfies for any $t\geq 0$ and $x\in H^{(1)}$
$$\int_{H^{(2)}}\chi(t,x,y)\mu^{x}(dy)=0.$$


Thanks to Lemma \ref{lemPoiss} above, we thus obtain the following solution to equation \eqref{defchi}
\begin{equation}\label{expressu1}
u_1(t,x,y)=\int_{0}^{+\infty}\E[\chi(t,x,Y_{x}(s,y))]ds.
\end{equation}
Moreover we are able to show regularity of $u_1$ with respect to $t$ and $x,y$.

The remainder $v^\epsilon=u^\epsilon-u_0-\epsilon u_1$ satisfies
\begin{equation}\label{estimdérivméthode}
(\partial_t-\frac{1}{\epsilon}L_1-L_2)v^\epsilon=\epsilon(L_2u_1-\frac{\partial u_1}{\partial t}).
\end{equation}

Due to non-integrability in $0$ of some bounds below, we introduce a parameter $\rho(\epsilon)=\epsilon^{1/\theta}\leq \epsilon$ (since $0<\theta\leq 1$); it satisfies $\rho(\epsilon)\rightarrow 0$ when $\epsilon\rightarrow 0$.

Using a variation of constant formula, we obtain
\begin{align*}
v^\epsilon(T,x,y)&=\E[v^\epsilon(\rho(\epsilon),X^\epsilon(T-\rho(\epsilon),x,y),Y^\epsilon(T-\rho(\epsilon),x,y))]\\
&+\epsilon\E[\int_{\rho(\epsilon)}^{T}(L_2u_1-\frac{\partial u_1}{\partial t})(t,X^\epsilon(T-t,x,y),Y^\epsilon(T-t,x,y))dt]
\end{align*}

By \eqref{expuepsilon}, and since $u_0=\overline{u}$, we then have
\begin{equation}\label{estimmé1}
\begin{aligned}
u^\epsilon(T,x,y)-\overline{u}(T,x,y)&=\epsilon u_1(T,x,y)\\
&+\E[v^\epsilon(\rho(\epsilon),X^\epsilon(T-\rho(\epsilon),x,y),Y^\epsilon(T-\rho(\epsilon),x,y))]\\
&+\epsilon\E[\int_{\rho(\epsilon)}^{T}(L_2u_1-\frac{\partial u_1}{\partial t})(t,X^\epsilon(T-t,x,y),Y^\epsilon(T-t,x,y))dt].
\end{aligned}
\end{equation}

The following estimates are proved below:
\begin{lemme}\label{lemu1final}
There exists a constant $C$ such that for any $0<t\leq T$, $x,y\in H$,
\begin{gather*}
|u_1(t,x,y)|\leq C(1+|x|+|y|)\\
|\frac{\partial u_1}{\partial t}(t,x,y)|\leq C(1+\frac{1}{t})(1+|x|+|y|)^2\\
|L_2u_1(t,x,y)|\leq C(1+|x|+|y|)(1+|Ax|).
\end{gather*}
\end{lemme}

Using estimates on $X^\epsilon$ and $Y^\epsilon$ proved in the appendix (see Propositions \ref{aprioribounds} and \ref{estimAX}) the first and the last expressions of \eqref{estimmé1} are bounded by
$$C\epsilon(1+|x|+|y|)+C\epsilon^{1-r/2}(1+|\log(\rho(\epsilon))|)(1+|x|_{(-A)^\theta}+|y|)^2,$$
which is dominated by $C\epsilon^{1-r}$ with the choice of $\rho(\epsilon)$ given above.

We notice that the Assumption $\theta>0$ is essential to control the part involving $|Ax|$.

We now explain how the central term of \eqref{estimmé1} is controlled; for that we estimate for any $x,y\in H$
\begin{equation}\label{estimmé2}
\begin{aligned}
v^\epsilon(\rho(\epsilon),x,y)&=u^\epsilon(\rho(\epsilon),x,y)-u_0(\rho(\epsilon),x)-\epsilon u_1(\rho(\epsilon),x,y)\\
&=-\epsilon u_1(\rho(\epsilon),x,y)\\
&+[u^\epsilon(\rho(\epsilon),x,y)-u^\epsilon(0,x,y)]-[u_0(\rho(\epsilon),x)-u_0(0,x)],
\end{aligned}
\end{equation}

since the initial condition $\phi$ is the same for $u^\epsilon$ and $\overline{u}$.

Using Lemma \ref{lemu1final}, the first term above is easily controlled by $C\epsilon(1+|x|+|y|)$. We now use another method to control the two other terms.

First, we use the definition \eqref{defu} of $\overline{u}=u_0$ to write
\begin{align*}
|u_0(\rho(\epsilon),x)-u_0(0,x)|&=|\int_{0}^{\rho(\epsilon)}\frac{\partial}{\partial t}u_0(t,x)dt|\\
&=|\int_{0}^{\rho(\epsilon)}\frac{\partial}{\partial t}\phi(\overline{X}(t,x))dt|\\
&=|\int_{0}^{\rho(\epsilon)}D\phi(\overline{X}(t,x)).\frac{d}{dt}\overline{X}(t,x)dt|\\
&\leq C\int_{0}^{\rho(\epsilon)}|\frac{d}{dt}\overline{X}(t,x)|dt.
\end{align*}

By definition of $\overline{X}$ (see \eqref{eqmoy}), and using Proposition \ref{Xbar3}, we get for any $t>0$
\begin{align*}
|\frac{d}{dt}\overline{X}(t,x)|\leq C(1+t^{\theta-1})(1+|x|_{(-A)^\theta}).
\end{align*}

As a consequence, since $\theta>0$, we get
$$|u_0(\rho(\epsilon),x)-u_0(0,x)|\leq C(\rho(\epsilon)+\frac{\rho(\epsilon)^{\theta}}{\theta})(1+|x|_{(-A)^\theta}).$$

The other expression is controlled in the same way; it is important to notice that the assumption that $\phi$ only depends on the slow variable $x$ is fundamental in this estimate.
\begin{align*}
|u^\epsilon(\rho(\epsilon),x,y)-u^\epsilon(0,x,y)|&=|\int_{0}^{\rho(\epsilon)}\frac{\partial}{\partial t}u^\epsilon(t,x,y)dt|\\
&=|\int_{0}^{\rho(\epsilon)}\frac{\partial}{\partial t}\E[\phi(X^\epsilon(t,x,y))]dt|\\
&=|\int_{0}^{\rho(\epsilon)}\E[D\phi(X^\epsilon(t,x,y)).\frac{d}{dt}X^\epsilon(t,x,y)]dt|\\
&\leq C\int_{0}^{\rho(\epsilon)}(\E|AX^\epsilon(t,x,y)|+1)dt.
\end{align*}

We now use the estimate on $\E|AX^\epsilon(t,x,y)|$ of Proposition \ref{estimAX} in the appendix, and we obtain
$$|u^\epsilon(\rho(\epsilon),x,y)-u^\epsilon(0,x,y)|\leq C(\rho(\epsilon)+\frac{\rho(\epsilon)^{\theta}}{\theta}+\epsilon^{-r/2}\rho(\epsilon))(1+|x|_{(-A)^\theta}+|y|).$$

Then by \eqref{estimmé2}, and using $\rho(\epsilon)=\epsilon^{1/\theta}\leq \epsilon$, we get
$$|v^\epsilon(\rho(\epsilon),x,y)|\leq C\epsilon^{1-r/2}(1+|x|_{(-A)^\theta}+|y|).$$

Then thanks to \eqref{estimmé1} and to Proposition \ref{coolestim}, we get for $\epsilon\leq 1$
$$|u^\epsilon(T,x,y)-\overline{u}(T,x,y)|\leq C\epsilon^{1-r}.$$

As expained at the end of Section \ref{sectionreduc}, we have indeed proved a bound on \eqref{weakN}. It is now enough to notice that the above constant $C$ is independent of dimension $N$, and to let $N$ go to $+\infty$, and Theorem \ref{weak} follows.

\section{Proof of Lemma \ref{lemu1final}}\label{proffestimsect}
We use results gathered in the appendix \ref{appendixXeYe} and \ref{appendixXbar}.

\subsection{Estimate of $u_1$}\label{sec u_1}

Since $u_1$ is defined by \eqref{expressu1}, by using Lemma \ref{lemPoiss} we have
$$|u_1(t,x,y)|\leq C(1+|y|^2)\|y\mapsto \chi(t,x,y)\|_{\infty}.$$
According to \eqref{defchi}, we indeed have for any $y\in H^{(2)}$
$$\chi(t,x,y)=<F(x,y)-\overline{F}(x),D_{x}u_0(t,x)>,$$ and therefore we just have to bound $|D_xu_0(t,x)|$, thanks to the following lemma:
\begin{lemme}\label{lem1}
For any $T\in]0,+\infty[$, there exists $C_0>0$ such that for any $0\leq t\leq T$ and $x\in H^{(1)}$
$$|D_xu_0(t,x)|_H\leq C_T\sup_{z\in H}|D\phi(z)|_H.$$
\end{lemme}

\underline{Proof}
$u_0$ is the solution of the equation
\begin{equation}\label{defu0}
\begin{gathered}
\frac{\partial u_0}{\partial t}(t,x)=<Ax+\overline{F}(x),D_xu_0(t,x)>\\
u_0(0,x)=\phi(x).
\end{gathered}
\end{equation}

We have a representation formula $u_0(t,x)=\phi(\overline{X}(t,x))$, where $\overline{X}$ is solution of (\ref{eqmoy}).

We can differentiate (\ref{eqmoy}) with respect to the initial condition $x$, and we have for any $h\in H^{(1)}$ $$D_xu_0(t,x).h=D\phi(\overline{X}(t,x)).\eta^{h}(t,x),$$ where $\eta^h(t,x)$ is the derivative of $\overline{X}$ with respect to $x$ in direction $h$, and is solution of the variational equation (\ref{defeta}) in the appendix.

Using Proposition \ref{eta1}, we get $|D_xu_0(t,x).h|\leq C_{T}\sup_{z\in H}|D\phi(z)||h|$, and taking the supremum over $h$ gives the result.
\begin{flushright}
$\Box$
\end{flushright}

Therefore, we obtain the first estimate of Lemma \ref{lemu1final}.

\subsection{Estimate of $\frac{\partial u_1}{\partial t}$.}\label{contdt}

First we check that
$$
\int_{H^{(2)}}\frac{\partial\chi}{\partial t}(t,x,z)\mu^x(dz)=\frac{\partial}{\partial t}\int_{H^{(2)}}\chi(t,x,z)\mu^x(dz)=0.
$$
By definition \eqref{expressu1} of $u_1$, it is easy to show that we can differentiate with respect to $t$, and that
\begin{equation}\label{expressdu1dt}
\frac{\partial u_1}{\partial t}(t,x,y)=\int_{0}^{+\infty}\E[\frac{\partial\chi}{\partial t}(t,x,Y_{x}(s,y))]ds.
\end{equation}
We then obtain
$$|\frac{\partial u_1}{\partial t}(t,x,y)|\leq C(1+|y|^2)\|y\mapsto \frac{\partial \chi}{\partial t}(t,x,y)\|_{\infty}.$$
Since by \eqref{defchi} we have
$$\frac{\partial\chi}{\partial t}(t,x,y)=<F(x,y)-\overline{F}(x),\frac{\partial }{\partial t}D_xu_0(t,x)>,$$
we just need to control $|\frac{\partial }{\partial t}D_xu_0(t,x)|$:

\begin{lemme}\label{lem3}
For any $T>0$, there exists $C_{T}>0$ such that for any $0<t\leq T$, $x\in H^{(1)}$ and $h\in H^{(1)}$ we have
$$|\frac{\partial }{\partial t}D_{x}u_0(t,x).h|\leq C(1+|x|_{H})(1+t^{-1})|h|.$$
\end{lemme}

\underline{Proof}
For any $h\in H^{(1)}$, we have
\begin{align*}
\frac{\partial }{\partial t}(D_xu_0(t,x).h)&=D^2\phi(\overline{X}(t,x))(\eta^h(t,x),\frac{d}{dt}\overline{X}(t,x))\\
&+D\phi(\overline{X}(t,x)).\frac{d}{dt}\eta^h(t,x).
\end{align*}

\begin{enumerate}
\item Thanks to Proposition \ref{Xbar1}, we have $|\eta^h(t,x)|\leq C(1+|x|)$ for any $t\geq 0$.

Moreover $\frac{d}{dt}\overline{X}(t,x)=A\overline{X}(t,x)+\overline{F}(\overline{X}(t,x))$.

On the one hand, $\overline{F}$ is bounded; on the other hand, thanks to Proposition \ref{Xbar3} we have
$$|A\overline{X}(t,x)|_{H}\leq C_{\theta}(1+t^{-1})(1+|x|_{H}).$$
Therefore
$$|\frac{d}{dt}\overline{X}(t,x)|_{H}\leq C(1+t^{-1})(1+|x|_{H}).$$

\item It remains to control
$$|\frac{d}{dt}\eta^h(t,x)|=|A\eta^h(t,x)+D\overline{F}(\overline{X}(t,x)).\eta^h(t,x)|.$$

Since $\overline{F}$ is Lipschitz continuous, and using Proposition \ref{Xbar1}, we get an estimate of the second term.

Moreover Proposition \ref{eta3} gives
$$|A\eta^h(t,x)|\leq C(t^{-1}+1)(1+|x|_{H})|h|.$$

Therefore $$|\frac{d}{dt}\eta^h(t,x)|\leq C(t^{-1}+1)(1+|x|_{H})|h|.$$

\item We then have for any $h\in H^{(1)}$
$$|\frac{\partial}{\partial t}(D_xu_0(t,x).h)|\leq C(t^{-1}+1)(1+|x|_{H})|h|.$$

\end{enumerate}
\begin{flushright}
$\Box$
\end{flushright}

We then obtain the second estimate of Lemma \ref{lemu1final}.

\subsection{Estimate of $L_2u_1$.}

To prove Lemma \ref{lemu1final}, it remains to control the part involving $L_2u_1$.

By definition of $L_2$, we have
\begin{equation}\label{expressL2u1}
L_2u_1(t,x,y)=<Ax+F(x,y),D_xu_1(t,x,y)>.
\end{equation}

Therefore we have to estimate $|D_xu_1(t,x,y)|$. We explain how $D_xu_1(t,x,y).h$ can be calculated for any $h\in H^{(1)}$.

Recall that $u_1$ defined by (\ref{expressu1}) satisfies
$$L_1u_1(t,x,y)=<\overline{F}(x)-F(x,y),D_{x}u_0(t,x)>=-\chi(t,x,y),$$
where we explicitely write the dependence of the operator $L_1$ in the two variables $x$ and $y$.

We fix $t\geq 0$, $x\in H^{(1)}$, $y\in H^{(2)}$, and $h\in H^{(1)}$. Then for any $\xi\neq 0$ we have
\begin{align*}
L_1(x,y)\frac{u_1(t,x+\xi h,y)-u_1(t,x,y)}{\xi}&=-\frac{\chi(t,x+\xi h,y)-\chi(t,s,y)}{\xi}\\
&-<\frac{G(x+\xi h,y)-G(x,y)}{\xi},D_yu_1(t,x+\xi h,y)>\\
&=:-\Gamma(t,x,y,h,\xi),
\end{align*}

where $\Gamma$ is regular with respect to $y$; therefore by using Lemma \ref{lemPoiss} we get
\begin{align*}
\frac{u_1(t,x+\xi h,y)-u_1(t,x,y)}{\xi}&-\int_{H^{(2)}}\frac{u_1(t,x+\xi h,y)-u_1(t,x,y)}{\xi}\mu^{x}(dy)\\
&=\int_{0}^{+\infty}\E[\Gamma(t,x,Y_x(s,y),h,\xi)]ds.
\end{align*}

We want to take the limit when $\xi\rightarrow 0$, in order to prove that we can differentiate, and to obtain an expression that we are able to control.

First, we notice that for any $t,x,y$ we have $\int_{H^{(2)}}u_1(t,x,y)\mu^{x}(dy)=0$; so we can write that
\begin{align*}
\int_{H^{(2)}}\frac{u_1(t,x+\xi h,y)-u_1(t,x,y)}{\xi}\mu^{x}(dy)&=-\int_{H^{(2)}}u_1(t,x+\xi h,y)\frac{V(x+\xi h,y)-V(x,y)}{\xi}\nu(dy),
\end{align*}

where $V(x,y):=\frac{1}{Z(x)}e^{2U(x,y)}$ (so that we have $\mu^x(dy)=V(x,y)\nu(dy)$).

When $\xi\rightarrow 0$, we obtain
\begin{align*}
\int_{H^{(2)}}\frac{u_1(t,x+\xi h,y)-u_1(t,x,y)}{\xi}\mu^{x}(dy)&\rightarrow \int_{H^{(2)}}u_1(t,x,y)D_xV(x,y).h\nu(dy)\\
&=\int_{H^{(2)}}u_1(t,x,y)H(x,y).hV(x,y)\nu(dy),
\end{align*}
where $H(x,y)=2D_xU(x,y)-2\int_{H^{(2)}}D_xU(x,z)\mu^x(dz)$.

Moreover $|\int_{H^{(2)}}u_1(t,x,y)H(x,y).hV(x,y)\nu(dy)|\leq C(1+|x|)|h|$.

Second, we look at the part involving $\Gamma$: we notice that when $\xi\rightarrow 0$,
$$\Gamma(t,x,y,h,\xi)\rightarrow \Theta(t,x,y).h,$$
where
\begin{equation}\label{expresstheta}
\Theta(t,x,y).h=D_x\chi(t,x,y).h+<D_xG(x,y).h,D_yu_1(t,x,y)>.
\end{equation}

Below, we prove the following estimate on the function $\Theta$:
\begin{lemme}\label{LemmaTheta}
There exists a constant $C$ such that for any $x\in H^{(1)}$, $t\geq 0$, $h\in H^{(1)}$ we have for any $y\in H^{(2)}$
$$|\Theta(t,x,y).h|\leq C(1+|y|^2)|h|.$$
\end{lemme}

We notice that for any $t,x,\xi,h$ we have by definition of $\Gamma$ $\int_{H^{(2)}}\Gamma(t,x,y,\xi,h)\mu^x(dy)=0$; then using the bound of the previous Lemma and the dominated convergence Theorem we obtain $\int_{H^{(2)}}\Theta(t,x,y).h\mu^x(dy)=0$ for any $x\in H^{(1)}$, $t\geq 0$, $h\in H^{(1)}$.
Using this result, Proposition \ref{propoexpy1y2} - with integration with respect to $\mu^x(dy_2)$ - and the estimate in the previous Lemma, we then see that $u_1$ can be differentiated with respect to $x$, and that the following formula holds:
\begin{equation}\label{expressdxu1}
D_xu_1(t,x,y).h=\int_{H^{(2)}}u_1(t,x,y)H(x,y).hV(x,y)\nu(dy)+\int_{0}^{+\infty}\E[\Theta(t,x,Y_x(s,y)).h]ds;
\end{equation}

According to Lemma \ref{LemmaTheta}, we do not know whether $\Theta$ is a bounded function, but we only know that it has quadratic growth. However, the result of Proposition \ref{propoexpy1y2} can easily be extended to such function.

Now we obtain that
$$|D_xu_1(t,x,y).h|\leq C(1+|y|^2)|h|$$
and therefore (see \ref{expressL2u1})
$$|L_2u_1(t,x,y)|\leq C(1+|y|^2)(1+|Ax|),$$
which is the third estimate of Lemma \ref{lemu1final}.

It remains to prove Lemma \ref{LemmaTheta}.


We fix $t\geq 0$, $x\in H^{(1)}$, $h\in H^{(1)}$, and $y,y'\in H^{(2)}$.

\begin{itemize}
\item On the one hand, $\chi$ being defined by (\ref{defchi}), we have
$$D_x\chi(t,x,y).h=<D_xF(x,y).h,D_{x}u_0(t,x)>+D_{xx}^2u_0(t,x).(h,F(x,y)).$$

Using the boundedness of the first derivative of $F$, and Lemma \ref{lem1}, we easily have
$$|<D_xF(x,y).h,D_{x}u_0(t,x)>|\leq C|h|.$$

The other part can be controlled thanks to the following Lemma:
\begin{lemme}
For any $0\leq t\leq T$, $x\in H^{(1)}$, $h,k\in H^{(1)}$, we have
$$|D_{xx}^{2}u_0(t,x).(h,k)|\leq C(T,\phi)|h|_{H}|k|_{H}.$$
\end{lemme}

\underline{Proof}
We have
\begin{gather*}
u_0(t,x)=\phi(\overline{X}(t,x))\\
D_xu_0(t,x).h=D\phi(\overline{X}(t,x)).(D_x\overline{X}(t,x).h)\\
\end{gather*}
and
\begin{align*}
D_{xx}^{2}u_0(t,x).(h,k)&=D^2\phi(\overline{X}(t,x))(D_x\overline{X}(t,x).h,D_x\overline{X}(t,x).k)\\
&+D\phi(\overline{X}(t,x)).(D_{xx}^{2}\overline{X}(t,x).(h,k)).
\end{align*}

Using Proposition \ref{eta1}, we control $\eta^{h}(t,x)=D_x\overline{X}(t,x).h$; moreover we notice that the second derivative $\xi^{h,k}(t,x):=D_{xx}^{2}\overline{X}(t,x).(h,k)$ satisifies equation (\ref{defxihk}); using Proposition \ref{propoxihk}, we get the result.
\begin{flushright}
$\Box$
\end{flushright}

Therefore $|D_x\chi(t,x,y).h|\leq C|h|$.

\item On the other hand,
$$|<D_xG(x,y).h,D_y u_1(t,x,y)>|\leq C|h||D_y u_1(t,x,y)|,$$

But we have proved in Lemma \ref{lemPoiss} how to control the derivatives of $u_1$ with respect to $y$: we obtain 
$$|D_yu_1(t,x,y).h|\leq C(1+|y|^2)|h|.$$

Therefore we have
$$|<D_xG(x,y).h,D_y u_1(t,x,y)>|\leq C(1+|y|^2)|h|,$$

and now the result is easily obtained:
$$|\Theta(t,x,y).h|\leq C(1+|y|^2)|h|.$$

\end{itemize}

\appendix

\section{Properties of $(X^\epsilon,Y^\epsilon)$}\label{appendixXeYe}


The results of this section only require Assumptions \ref{hypAB}, \ref{hypF} and \ref{hypG}; in particular no dissipativity is assumed.

The first important property is the control of moments of any order:

\begin{propo}\label{aprioribounds}
For any $1\leq p<+\infty$, there exists $c_{p}>0$ such that for any $(x,y)\in H^2$, $t\geq 0$ and $\epsilon>0$
$$\E[|X^{\epsilon}(t)|_{H}^{p}]\leq c_{p}(1+e^{-\lambda t}|x|_{H}^{p}) \quad \text{and} \quad \E[|Y^{\epsilon}(t)|_{H}^{p}]\leq c_{p}(1+e^{-\mu t}|y|_{H}^{p}).$$

\end{propo}

We can also give bounds on the moments with respect to $|\hspace{3pt}.\hspace{3pt}|_{(-A)^a}$ and $|\hspace{3pt}.\hspace{3pt}|_{(-B)^b}$ norms, for $0<a<1$ and $0<b<1/4$ (the case $a=1$ is treated in Proposition \ref{estimAX} below).

\begin{propo}\label{coolestim}
For any $p\geq 1$, $a\in(0,1)$, $b\in (0,1/4)$, there exists $C_{p,a,b}>0$ such that for any $x\in D(-A)^a$ and $y\in D(-B)^b$, we have:
$$\E|X^\epsilon(t,x,y)|_{(-A)^a}^p\leq C_p(1+|x|_{(-A)^a}^{p}) \quad \text{and} \quad \E|Y^\epsilon(t,x,y)|_{(-B)^b}^p\leq C_p(1+|y|_{(-B)^b}^{p}).$$
\end{propo}

We now give some regularity estimates of $X^\epsilon$ and $Y^\epsilon$ in the time variable. We do not assume any regularity assumption on $x$ or $y$; as a consequence, we obtain singularities at the origin, which are integrable.

\begin{propo}\label{reg}
For any $0<r<1$, there exists $C_r>0$ such that for any $x,y\in H$, for any $0<s\leq t$ and $\epsilon>0$ we have
$$\left(\E|X^{\epsilon}(t)-X^{\epsilon}(s)|_{H}^{2}\right)^{1/2}\leq C_r|t-s|^{1-r}(1+\frac{1}{s^{1-r}})(1+|x|_{H}).$$
\end{propo}

\underline{Proof}
If we fix $0<s\leq t$, $x,y\in H$, we have
\begin{align*}
X^{\epsilon}(t)-X^{\epsilon}(s)&=e^{tA}x-e^{sA}x\\
&+\int_{0}^{t}e^{(t-\sigma)A}F(X^{\epsilon}(\sigma),Y^{\epsilon}(\sigma))d\sigma-\int_{0}^{s}e^{(s-\sigma)A}F(X^{\epsilon}(\sigma),Y^{\epsilon}(\sigma))d\sigma.
\end{align*}

For the first term, if $x=\sum_{k=0}^{+\infty}x_ke_k$, we can use Proposition \ref{proporegul} to get
$$|e^{tA}x-e^{sA}x|_{H}\leq C_{r}\frac{(t-s)^{1-r}}{s^{1-r}}|x|_{H}.$$

For the second term, we use the following decomposition:
\begin{align*}
\int_{0}^{t}e^{(t-\sigma)A}F(X^{\epsilon}(\sigma),Y^{\epsilon}(\sigma))d\sigma-&\int_{0}^{s}e^{(s-\sigma)A}F(X^{\epsilon}(\sigma),Y^{\epsilon}(\sigma))d\sigma\\
&=\int_{s}^{t}e^{(t-\sigma)A}F(X^{\epsilon}(\sigma),Y^{\epsilon}(\sigma))d\sigma\\
&+\int_{0}^{s}(e^{(t-\sigma)A}-e^{(s-\sigma)A})F(X^{\epsilon}(\sigma),Y^{\epsilon}(\sigma))d\sigma.
\end{align*}

First, by the Cauchy-Schwarz inequality, we have
\begin{align*}
\E|\int_{s}^{t}e^{(t-\sigma)A}F(X^{\epsilon}(\sigma),Y^{\epsilon}(\sigma))d\sigma|_{H}^{2}&\leq
(t-s)\E\int_{s}^{t}|e^{(t-\sigma)A}F(X^{\epsilon}(\sigma),Y^{\epsilon}(\sigma))|_{H}d\sigma\\
&\leq C(t-s)^2,
\end{align*}
since $F$ is assumed to be bounded.

Second, we use the second inequality of Proposition \ref{proporegul} to control the last expression:
\begin{align*}
\E|\int_{0}^{s}e^{(s-\sigma)A}&(e^{(t-s)A}-I)F(X^{\epsilon}(\sigma),Y^{\epsilon}(\sigma))d\sigma|_{H}^{2}\\
&\leq \E[\int_{0}^{s}|\left(e^{(t-\sigma)A}-e^{(s-\sigma)A}\right)F(X^{\epsilon}(\sigma),Y^{\epsilon}(\sigma))d\sigma|_{H}]^2\\
&\leq C_{r}^{2}(t-s)^{2(1-r)}\E(\int_{0}^{s}\frac{e^{-\frac{\lambda}{2}(s-\sigma)}}{(s-\sigma)^{1-r}}|F(X^{\epsilon}(\sigma),Y^{\epsilon}(\sigma))|_{H}d\sigma)^2\\
&\leq C_{r}^{2}(t-s)^{2(1-r)},
\end{align*}
since $\int_{0}^{+\infty}\frac{e^{-\frac{\lambda}{2} s}}{s^{1-r}}ds<+\infty$.
\begin{flushright}
$\Box$
\end{flushright}

\begin{propo}\label{regY}
For any $0<r<1/4$, there exists a constant $C_r$ such that if $x,y\in H$, then for any $0<s<t$ and $\epsilon>0$
$$\E|Y^\epsilon(t)-Y^\epsilon(s)|^2\leq C(1+|x|_{H}^{2}+|y|_{H}^{2})[\left(\frac{t-s}{s}\right)^{2r}+\left(\frac{t-s}{\epsilon}\right)^{2r}].$$
\end{propo}

\underline{Proof}

\begin{itemize}
\item For any $0<s<t$,
\begin{align*}
Y^\epsilon(t)-Y^\epsilon(s)&=(e^{\frac{t}{\epsilon}B}-e^{\frac{s}{\epsilon}B})y\\
&+\frac{1}{\epsilon}\int_{s}^{t}e^{\frac{(t-\sigma)}{\epsilon}B}G(X^\epsilon(\sigma),Y^\epsilon(\sigma))d\sigma\\
&+\frac{1}{\epsilon}\int_{0}^{s}\left(e^{\frac{(t-\sigma)}{\epsilon}B}-e^{\frac{(s-\sigma)}{\epsilon}B}\right)G(X^\epsilon(\sigma),Y^\epsilon(\sigma))d\sigma\\
&+W^{\epsilon,B}(t)-W^{\epsilon,B}(s),
\end{align*}

where $W^{\epsilon,B}(r)=\frac{1}{\sqrt{\epsilon}}\int_{0}^{r}e^{\frac{(r-\sigma)}{\epsilon}B}dW(\sigma)$. We remark that only the last expression can not bounded almost surely (since we assume that $G$ is bounded).

\item For the first term, using the second inequality of Proposition \ref{proporegul}, we have for any $0<s<t$
$$|(e^{\frac{t}{\epsilon}B}-e^{\frac{s}{\epsilon}B})y|_{H}\leq C_r(\frac{t-s}{s})^{r}|y|_{H}.$$

\item For the second term, we have for any $0<s<t$
\begin{align*}
|\frac{1}{\epsilon}\int_{s}^{t}e^{\frac{(t-\sigma)}{\epsilon}B}G(X^\epsilon(\sigma),Y^\epsilon(\sigma))d\sigma|&\leq \frac{1}{\epsilon}\int_{s}^{t}|e^{\frac{(t-\sigma)}{\epsilon}B}|_{\mathcal{L}(H)}\|G\|_{\infty}d\sigma\\
&\leq \frac{C}{\epsilon}\int_{s}^{t}e^{-\mu(t-\sigma)/\epsilon}d\sigma\\
&\leq C\int_{0}^{(t-s)/\epsilon}e^{-\mu\sigma}d\sigma\\
&\leq C \frac{(t-s)^r}{\epsilon^r}.
\end{align*}

\item For the third term, we use the second estimate of Proposition \ref{proporegul}, and we have for any $0<s<t$
\begin{align*}
|\frac{1}{\epsilon}\int_{0}^{s}\left(e^{\frac{(t-\sigma)}{\epsilon}B}-e^{\frac{(s-\sigma)}{\epsilon}B}\right)G(X^\epsilon(\sigma),Y^\epsilon(\sigma))d\sigma|&\leq \frac{1}{\epsilon}\int_{0}^{s}|e^{\frac{(t-\sigma)}{\epsilon}B}-e^{\frac{(s-\sigma)}{\epsilon}B}|_{\mathcal{L}(H)}\|G\|_{\infty}d\sigma\\
&\leq \frac{C_r}{\epsilon}\int_{0}^{s}\frac{(t-s)^r}{(s-\sigma)^r}e^{-\frac{\mu(s-\sigma)}{2\epsilon}}d\sigma\\
&\leq C_r\frac{(t-s)^{r}}{\epsilon^{r}}\int_{0}^{+\infty}\frac{1}{\sigma^r}e^{-\frac{\mu\sigma}{2}}d\sigma.
\end{align*}

\item For the fourth term, we have for any $0<s<t$,
\begin{align*}
\E|W^{\epsilon,B}(t)-W^{\epsilon,B}(s)|^2&=\E|\frac{1}{\sqrt{\epsilon}}\int_{s}^{t}e^{(t-\sigma)B/\epsilon}dW(\sigma)+\frac{1}{\sqrt{\epsilon}}\int_{0}^{s}(e^{(t-\sigma)B/\epsilon}-e^{(s-\sigma)B/\epsilon})dW(\sigma)|^2\\
&=\E|\frac{1}{\sqrt{\epsilon}}\int_{s}^{t}e^{(t-\sigma)B/\epsilon}dW(\sigma)|^2\\
&+\E|\frac{1}{\sqrt{\epsilon}}\int_{0}^{s}(e^{(t-\sigma)B/\epsilon}-e^{(s-\sigma)B/\epsilon})dW(\sigma)|^2\\
&=\frac{1}{\epsilon}\int_{s}^{t}|e^{(t-\sigma)B/\epsilon}|_{\mathcal{L}_2(H)}^{2}d\sigma\\
&+\frac{1}{\epsilon}\int_{0}^{t}|e^{(t-\sigma)B/\epsilon}-e^{(s-\sigma)B/\epsilon}|_{\mathcal{L}_2(H)}^{2}d\sigma.
\end{align*}

On the one hand,
\begin{align*}
\frac{1}{\epsilon}\int_{s}^{t}|e^{(t-\sigma)B/\epsilon}|_{\mathcal{L}_2(H)}^{2}d\sigma&=\frac{1}{\epsilon}\int_{s}^{t}\sum_{k=0}^{+\infty}e^{-2(t-\sigma)\mu_k/\epsilon}d\sigma\\
&=\sum_{k=0}^{+\infty}\int_{0}^{(t-s)/\epsilon}e^{-2\sigma\mu_k}d\sigma\\
&=\sum_{k=0}^{+\infty}\frac{1}{2\mu_k}(1-e^{-2\mu_k(t-s)/\epsilon})\\
&\leq C_\zeta\sum_{k=0}^{+\infty}\frac{\mu_{k}^{2r}}{\mu_k}\left(\frac{t-s}{\epsilon}\right)^{2r},
\end{align*}

and we know (by Assumption \ref{hypAB}) that the above sum is finite if and only if $r<1/4$;

on the other hand,
\begin{align*}
\frac{1}{\epsilon}\int_{0}^{s}|e^{(t-\sigma)B/\epsilon}-e^{(s-\sigma)B/\epsilon}|_{\mathcal{L}_2(H)}^{2}d\sigma&=\frac{1}{\epsilon}\int_{0}^{s}\sum_{k=0}^{+\infty}e^{-2(s-\sigma)\mu_k/\epsilon}(1-e^{-(t-s)\mu_k/\epsilon})^2d\sigma\\
&\leq C_\zeta\sum_{k=0}^{+\infty}\left(\frac{t-s}{\epsilon}\right)^{2r}\frac{\mu_{k}^{2r}}{\mu_k}(1-e^{-2s\mu_k}).
\end{align*}

\end{itemize}
\begin{flushright}
$\Box$
\end{flushright}

Finally the following Proposition gives a control for $AX^\epsilon$. We assume $\theta>0$, even if the proof is valid for $\theta=0$.
\begin{propo}\label{estimAX}
For any $0<r<1$, there exists a constant $C_{r}$ such that if $x\in D((-A)^\theta)$ and $y\in H$, then for any $t>0$ and $\epsilon>0$
$$(\E[|AX^\epsilon(t)|^2])^{1/2}\leq C_r(1+t^{\theta-1})|x|_{(-A)^\theta}+C_{r}(1+\epsilon^{-\frac{r}{2}})(1+|x|_{H}+|y|_{H}).$$
\end{propo}

\underline{Proof}
We remark that in Lemma $4.4$ of \cite{Ce-F} $\E|AX^\epsilon(t)|$ is controlled, but the same approach gives the result for $\E|AX^\epsilon(t)|^2$.
We have
\begin{align*}
X^\epsilon(t)&=e^{tA}x+\int_{0}^{t}e^{(t-s)A}F(X^\epsilon(s),Y^\epsilon(s))ds\\
&=e^{tA}x+\int_{0}^{t}e^{(t-s)A}F(X^\epsilon(t),Y^\epsilon(t))ds\\
&+\int_{0}^{t}e^{(t-s)A}\left(F(X^\epsilon(s),Y^\epsilon(s))-F(X^\epsilon(t),Y^\epsilon(t))\right)ds.
\end{align*}

For the first term, we have for any $t>0$
$$|Ae^{tA}x|\leq Ct^{\theta-1}|x|_{(-A)^\theta}.$$

For the second term, we have
$$|A\int_{0}^{t}e^{(t-s)A}F(X^\epsilon(t),Y^\epsilon(t))ds|=|(e^{tA}-I)F(X^\epsilon(t),Y^\epsilon(t))|\leq C.$$

For the third term, we have
\begin{multline*}
|A\int_{0}^{t}e^{(t-s)A}\left(F(X^\epsilon(s),Y^\epsilon(s))-F(X^\epsilon(t),Y^\epsilon(t))\right)ds|\\
\leq \int_{0}^{t}\frac{Ce^{-\frac{\lambda}{2}(t-s)}}{t-s}\left( |X^\epsilon(s)-X^\epsilon(t)|+|Y^\epsilon(s)-Y^\epsilon(t)|\right)ds.
\end{multline*}

Using Minkowski inequality, we get
\begin{gather*}
\E\left(\int_{0}^{t}\frac{Ce^{-\frac{\lambda}{2}(t-s)}}{t-s}|X^\epsilon(s)-X^\epsilon(t)|ds\right)^2\leq \left(\int_{0}^{t}\frac{Ce^{-\frac{\lambda}{2}(t-s)}}{t-s}(\E|X^\epsilon(t)-X^\epsilon(s)|^2)^{1/2}ds\right)^2\\
\E\left(\int_{0}^{t}\frac{Ce^{-\frac{\lambda}{2}(t-s)}}{t-s}|Y^\epsilon(s)-Y^\epsilon(t)|ds\right)^2\leq \left(\int_{0}^{t}\frac{Ce^{-\frac{\lambda}{2}(t-s)}}{t-s}(\E|Y^\epsilon(t)-Y^\epsilon(s)|^2)^{1/2}ds\right)^2.
\end{gather*}

Using Propositions \ref{reg} and \ref{regY}, we obtain a regularity result which gives convergent integrals. It is then easy to conclude.
\begin{flushright}
$\Box$
\end{flushright}

\section{Properties of $\overline{X}$}\label{appendixXbar}

Again the results of this section only require Assumptions \ref{hypAB}, \ref{hypF} and \ref{hypG}; in particular no dissipativity is assumed.

Recall that $\overline{X}(t,x)$ is defined via \eqref{eqmoy}.

\begin{propo}\label{Xbar1}
There exists $C>0$ such that for any $x\in H$ and any $t\geq 0$
$$|\overline{X}(t,x)|\leq C(1+e^{-\lambda t}|x|).$$
\end{propo}

\underline{Proof}
We use the mild representation formula: for any $t\geq 0$ and $x\in H$,
\begin{align*}
|\overline{X}(t,x)|&=|e^{tA}x+\int_{0}^{t}e^{(t-s)A}\overline{F}(\overline{X}(s,x))ds|\\
&\leq e^{-\lambda t}x+\int_{0}^{t}e^{-\lambda(t-s)}|\overline{F}(\overline{X}(s,x))|ds\\
&\leq C(1+e^{-\lambda t}|x|),
\end{align*}
since $\overline{F}$ is bounded.
\begin{flushright}
$\Box$
\end{flushright}

\begin{propo}\label{Xbar2}
For any $0<r<1$ and $0<\theta\leq 1$, there exists $C_r>0$ such that for any $x\in H$, for any $0<s\leq t$, we have
$$|\overline{X}(t,x)-\overline{X}(s,x)|\leq C_r|t-s|^{1-r}(1+\frac{1}{s^{1-r}})(1+|x|_{H}).$$
\end{propo}

\underline{Proof}
\begin{itemize}
\item If $0\leq s<t\leq T$, we can write
\begin{align*}
\overline{X}(t,x)-\overline{X}(s,x)&=(e^{tA}-e^{sA})x\\
&+\int_{s}^{t}e^{(t-\sigma)A}\overline{F}(\overline{X}(\sigma,x))d\sigma\\
&+\int_{0}^{s}(e^{(t-\sigma)A}-e^{(s-\sigma)A})\overline{F}(\overline{X}(\sigma,x))d\sigma.
\end{align*}

\item For the first term, it is easy to see that $|(e^{tA}-e^{sA})x|_{H}\leq C_{r}|t-s|^{1-r}(1+\frac{1}{s^{1-r}})|x|_{H}$.

\item For the second term, since $\overline{F}$ is bounded we have $|\int_{s}^{t}e^{(t-\sigma)A}\overline{F}(\overline{X}(\sigma,x))d\sigma|_{H}\leq C(t-s)$.

\item For the third term, we have
\begin{align*}
|\int_{0}^{s}(e^{(t-\sigma)A}-e^{(s-\sigma)A})\overline{F}(\overline{X}(\sigma,x))d\sigma|_{H}&\leq C_{r}\int_{0}^{s}\frac{e^{-\frac{\lambda}{2}(s-\sigma)}}{(s-\sigma)^{1-r}}(t-s)^{1-r}|\overline{F}(\overline{X}(\sigma,x))|_{H}d\sigma\\
&\leq C_{r}(t-s)^{1-r},
\end{align*}
\end{itemize}

\begin{flushright}
$\Box$
\end{flushright}

\begin{propo}\label{Xbar3}
For any $0<\theta\leq 1$, there exists $C(\theta)>0$ such that if $x\in D(-A)^\theta$, then for any $t>0$
$$|A\overline{X}(t,x)|_{H}\leq C_{\theta}(1+t^{\theta-1})(1+|x|_{(-A)^\theta}).$$
\end{propo}

\underline{Proof}
We first write that for any $t\geq 0$
\begin{align*}
\overline{X}(t,x)&=e^{tA}x+\int_{0}^{t}e^{(t-s)A}\overline{F}(\overline{X}(s,x))ds\\
&=e^{tA}x+\int_{0}^{t}e^{(t-s)A}\overline{F}(\overline{X}(t,x))ds\\
&+\int_{0}^{t}e^{(t-s)A}(\overline{F}(\overline{X}(s,x))-\overline{F}(\overline{X}(t,x)))ds.
\end{align*}

We have $|Ae^{tA}x|_{H}\leq C|x|_{(-A)^\theta}t^{\theta-1}$.

For the second term, we have
\begin{align*}
|A\int_{0}^{t}e^{(t-s)A}\overline{F}(\overline{X}(t,x))ds|_{H}&=|(e^{tA}-I)\overline{F}(\overline{X}(t,x))|_{H}\\
&\leq |\overline{F}(\overline{X}(t,x))|_{H}\\
&\leq C.
\end{align*}

The third term can be controlled by
\begin{align*}
|A\int_{0}^{t}e^{(t-s)A}(\overline{F}(\overline{X}(s,x))-\overline{F}(\overline{X}(t,x)))ds|&\leq C\int_{0}^{t}\frac{e^{-c(t-s)}}{t-s}|\overline{X}(t,x)-\overline{X}(s,x)|_{H}ds.
\end{align*}

In order to get a convergent integral, we use the regularity result of $\overline{X}$ proved in Proposition \ref{Xbar2}; therefore we obtain the result.

\begin{flushright}
$\Box$
\end{flushright}

The next three Propositions deal with $\eta^{h}(t,x)$ the derivative of $\overline{X}(t,x)$ with respect to $x$ in direction $h\in H$, at time $t$: it is the solution of
\begin{equation}\label{defeta}
\begin{gathered}
\frac{d\eta^h(t,x)}{dt}=A\eta^h(t,x)+D\overline{F}(\overline{X}(t,x)).\eta^h(t,x)\\
\eta^h(0,x)=h.
\end{gathered}
\end{equation}

Notice that we have to consider a finite horizon $T>0$.

\begin{propo}\label{eta1}
For any $T>0$, there exists $C_T>0$ such that for any $x\in H$, $h\in H$ and $0< t\leq T$
\begin{gather*}
|\eta^h(t,x)|\leq C_T|h|\\
|\eta^h(t,x)|_{(-A)^\eta}\leq C_T(1+\frac{1}{t^\eta})|h|.
\end{gather*}

\end{propo}

\underline{Proof}
We use that $A$ is a negative operator to prove say that for any $t\geq 0$
\begin{align*}
\frac{1}{2}\frac{d |\eta^h(t,x)|^2}{dt}&=<A\eta^h(t,x),\eta^h(t,x)>+<D\overline{F}(\overline{X}(t,x)).\eta^h(t,x),\eta^h(t,x)>\\
&\leq [\overline{F}]_{\text{Lip}}|\eta^h(t,x)|^2\\
&\leq C|\eta^h(t,x)|^2.
\end{align*}

Gronwall Lemma then yields the first estimate. The second one is proved by using the mild formulation for $\eta^h(t,x)$:
$$\eta^h(t,x)=e^{tA}h+\int_{0}^{t}e^{(t-s)A}D\overline{F}(\overline{X}(s,x)).\eta^h(s,x)ds;$$
thanks to the previous estimate, the integral is bounded by a constant, while $|e^{tA}h|_{(-A)^\eta}\leq \frac{C}{t^\eta}|h|_H$ (see Proposition \ref{proporegul}).

\begin{flushright}
$\Box$
\end{flushright}

\begin{propo}\label{eta2}
For any $T>0$, $0<r<1$, there exists $C_{T,r}>0$ such that for any $x\in H$, $h\in H$ and $0<s\leq t\leq T$
$$|\eta^h(t,x)-\eta^h(s,x)|\leq C_{T,r}(t-s)^{1-r}(1+\frac{1}{s^{1-r}})|h|.$$
\end{propo}

\underline{Proof}
\begin{itemize}
\item For $0\leq s<t\leq T$ we can write that
\begin{align*}
\eta^h(t,x)-\eta^h(s,x)&=(e^{tA}-e^{sA})h+\int_{s}^{t}e^{(t-\sigma)A}D\overline{F}(\overline{X}(\sigma,x)).\eta^h(\sigma,x)d\sigma\\
&+\int_{0}^{s}(e^{(t-\sigma)A}-e^{(s-\sigma)A})D\overline{F}(\overline{X}(\sigma,x)).\eta^h(\sigma,x)d\sigma.
\end{align*}

\item For the first term, we can see that $|(e^{tA}-e^{sA})h|\leq C_r\frac{(t-s)^{1-r}}{s^{1-r}}|h|$.

\item For the second term, we simply have
$$|\int_{s}^{t}e^{(t-\sigma)A}D\overline{F}(\overline{X}(\sigma,x)).\eta^h(\sigma,x)d\sigma|\leq C(t-s)|h|_{H}.$$

\item For the third term,
\begin{align*}
|\int_{0}^{s}(e^{(t-\sigma)A}-e^{(s-\sigma)A})&D\overline{F}(\overline{X}(\sigma,x)).\eta^h(\sigma,x)d\sigma|\\&\leq
C_{\delta}\int_{0}^{s}\frac{(t-s)^{1-r}}{(s-\sigma)^{1-r}}|D\overline{F}(\overline{X}(\sigma,x)).\eta^h(\sigma,x)|_{H}d\sigma\\
&\leq C_{r,T}(t-s)^{1-r}|h|_{H}.
\end{align*}

\end{itemize}
\begin{flushright}
$\Box$
\end{flushright}

\begin{propo}\label{eta3}
For any $T>0$, there exists $C_T$ such that for any $x\in H$, $h\in H$ and $0<t\leq T$
$$|A\eta^h(t,x)|\leq C_T(t^{-1}+1)(1+|x|)|h|.$$
\end{propo}

\underline{Proof}
For any $t\geq 0$,
\begin{align*}
\eta^h(t,x)&=e^{tA}h+\int_{0}^{t}e^{(t-s)A}D\overline{F}(\overline{X}(s,x)).\eta^h(s,x)ds\\
&=e^{tA}h+\int_{0}^{t}e^{(t-s)A}D\overline{F}(\overline{X}(t,x)).\eta^h(t,x)ds\\
&+\int_{0}^{t}e^{(t-s)A}(D\overline{F}(\overline{X}(s,x)).\eta^h(s,x)-D\overline{F}(\overline{X}(t,x)).\eta^h(t,x))ds.
\end{align*}

For the first term, we have $|Ae^{tA}h|_{H}\leq C t^{-1}|h|$.

For the second term,
\begin{align*}
|A\int_{0}^{t}e^{(t-s)A}D\overline{F}(\overline{X}(t,x)).\eta^h(t,x)ds|_{H}&=|(e^{tA}-I)D\overline{F}(\overline{X}(t,x)).\eta^h(t,x)|_{H}\\
&\leq 2|D\overline{F}(\overline{X}(t,x)).\eta^h(t,x)|_{H}\\
&\leq C|\eta^h(t,x)|_{H}\\
&\leq C|h|_{H}.
\end{align*}

For the third term, we have
\begin{multline*}
|A\int_{0}^{t}e^{(t-s)A}(D\overline{F}(\overline{X}(s,x)).\eta^h(s,x)-D\overline{F}(\overline{X}(t,x)).\eta^h(t,x))ds|_{H}\\
\leq \int_{0}^{t}\frac{C}{t-s}|D\overline{F}(\overline{X}(t,x)).\eta^h(t,x)-D\overline{F}(\overline{X}(s,x)).\eta^h(s,x)|_{H}ds.
\end{multline*}

To get a convergent integral, we need to show some regularity property.

For any $0\leq s<t\leq T$,

\begin{align*}
D\overline{F}(\overline{X}(t,x)).\eta^h(t,x)-D\overline{F}(\overline{X}(s,x)).\eta^h(s,x)&=[D\overline{F}(\overline{X}(t,x))-D\overline{F}(\overline{X}(s,x))].\eta^h(t,x)\\
&+D\overline{F}(\overline{X}(s,x)).(\eta^h(t,x)-\eta^h(s,x)).
\end{align*}

On the one hand, using Proposition \ref{FbarLip} on the regularity of $\overline{F}$, we have
\begin{align*}
|[D\overline{F}(\overline{X}(t,x))-D\overline{F}(\overline{X}(s,x))].\eta^h(t,x)|&\leq C|\overline{X}(t,x)-\overline{X}(s,x)||\eta^h(t,x)|_{(-A)^\eta}\\
&\leq C(1+|x|)(t-s)^{r}(1+\frac{1}{s^r})(1+\frac{1}{s^\eta})|h|_{H},
\end{align*}

thanks to Propositions \ref{Xbar2} and \ref{eta1}. Here $r$ must satisfy $r>0$ and $\eta+r<1$.

On the other hand, using Proposition \ref{eta2},
\begin{align*}
|D\overline{F}(\overline{X}(s,x)).(\eta^h(t,x)-\eta^h(s,x))|_{H}&\leq C|\eta^h(t,x)-\eta^h(s,x)|_{H}\\
&\leq C|h||t-s|^{1-r}(1+\frac{1}{s^{1-r}}).
\end{align*}

By integration, we then obtain the result.
\begin{flushright}
$\Box$
\end{flushright}

Finally we focus on $\xi^{h,k}(t,x)$ the second derivative of $\overline{X}(t,x)$ with respect to $x$ in directions $h,k\in H$, at time $t$: it is solution of
\begin{equation}\label{defxihk}
\begin{aligned}
\frac{d\xi^{h,k}(t,x)}{dt}&=A\xi^{h,k}(t,x)+D_x\overline{F}(\overline{X}(t,x)).(\xi^{h,k}(t,x))\\
&+D_{xx}^{2}\overline{F}(\overline{X}(t,x)).(\eta^h(t,x),\eta^k(t,x)).
\end{aligned}
\end{equation}

\begin{propo}\label{propoxihk}
For any $T>0$, there exists $C_T>0$ such that for any $x\in H$, $h,k\in H$ and $0\leq t\leq T$
$$|\xi^{h,k}(t,x)|\leq C_T|h||k|.$$
\end{propo}

\underline{Proof}
We have - since $A$ is negative, and using the estimates of Proposition \ref{FbarLip}:
\begin{align*}
\frac{1}{2}\frac{d|\xi^{h,k}(t,x)|^2}{dt}&\leq |D_x\overline{F}(\overline{X}(t,x))||\xi^{h,k}(t,x)|^2\\
&+C|\eta^h(t,x)||\eta^k(t,x)|_{(-A)^\eta}|\xi^{h,k}(t,x)|\\
&\leq C|\xi^{h,k}(t,x)|^2+C|\eta^h(t,x)|^2|\eta^k(t,x)|_{(-A)^\eta}^{2},
\end{align*}

Using Proposition \ref{eta1}, the Assumption $\eta<\frac{1}{2}$, and the Gronwall Lemma, we get the result.

\begin{flushright}
$\Box$
\end{flushright}

\section{Properties of the auxiliary function $\tilde{F}$}\label{sectionfunctionsaux}

The results of this section are used only for the proof of the strong convergence Theorem \ref{strong}. Here we need the strict dissipativity Assumption \ref{strictdiss}. 

In the proof of Lemma \ref{lemstrong2}, we need to use an auxiliary function $\tilde{F}$ - see definition \ref{defiauxfunctions}.

Thanks to Proposition \ref{convexp}, we get:
\begin{propo}\label{tildeFinf}
There exists $c>0$, $C>0$ such that for any $(x,y)\in H^2$ and $t\geq 0$,
$$|\tilde{F}(x,y,t)|\leq Ce^{-ct}(1+|x|_{H}+|y|_{H}).$$
\end{propo}

We also need the following estimate on the Lipschitz constant of $\tilde{F}$ with respect to $x$, which depends on the regularity assumptions made on $F$ and $G$ - see Assumptions \ref{hypF} and \ref{hypG}:

\begin{propo}\label{tildeFLip}
There exists $c>0$, $C>0$, such that for any $x_1,x_2,y\in H$ and $t\geq 0$
$$|\tilde{F}(x_1,y,t)-\tilde{F}(x_2,y,t)|\leq C(1+|y|)e^{-ct}(1+\frac{1}{t^\eta})|x_1-x_2|.$$
\end{propo}

\underline{Proof}
For any $t_0>0$, we define the following function:
$$\tilde{F}_{t_0}(x,y,t)=\hat{F}(x,y,t)-\hat{F}(x,y,t+t_0),$$

where $\hat{F}(x,y,t):=\E F(x,Y_{x}(t,y))$.

We claim that it satisfies the following properties:
\begin{itemize}
\item $\tilde{F}_{t_0}(x,y,t)\rightarrow \tilde{F}(x,y,t)$ when $t_0\rightarrow +\infty$.
\item For any $t_0$, for any $x,y,t$ and any $h$, $\tilde{F}_{t_0}$ is differentiable with respect to $x$ at $(x,y,t)$ and in direction $h\in H$.
\item We have $|D_x\tilde{F}_{t_0}(x,y,t).h|\leq Ce^{-ct}(1+\frac{1}{t^\eta})(1+|y|)|h|$, $C$ being independent of $t_0$.
\end{itemize}

The first two ones are obvious, thanks to regularity properties of $F$; moreover as soon as we have the third property, the proof of the Proposition can be finished as follows: if we fix $x_1,x_2,y,t,h$, then for any $t_0>0$
$$|\tilde{F}_{t_0}(x_1,y,t)-\tilde{F}_{t_0}(x_2,y,t)|\leq Ce^{-ct}(1+\frac{1}{t^\eta})(1+|y|)|x_1-x_2|.$$

Letting $t_0\rightarrow +\infty$, we get
$$|\tilde{F}(x_1,y,t)-\tilde{F}(x_2,y,t)|\leq Ce^{-ct}(1+\frac{1}{t^\eta})(1+|y|)|x_1-x_2|.$$

It remains to estimate $|D_x\tilde{F}_{t_0}(x,y,t).h|$ for any $h\in H$.

First we notice that thanks to the Markov property we have
\begin{align*}
\tilde{F}_{t_0}(x,y,t)&=\hat{F}(x,y,t)-\hat{F}(x,y,t+t_0)\\
&=\hat{F}(x,y,t)-\E F(x,Y_{x}(t+t_0,y))\\
&=\hat{F}(x,y,t)-\E\hat{F}(x,Y_{x}(t_0,y),t).
\end{align*}

Therefore we have for any $h$
\begin{align*}
D_x\tilde{F}_{t_0}(x,y,t).h&=D_x\hat{F}(x,y,t).h-\E D_x\left(\hat{F}(x,Y_{x}(t_0,y),t)\right).h\\
&=D_x\hat{F}(x,y,t).h-\E D_x\hat{F}(x,Y_{x}(t_0,y),t).h-\E D_y\hat{F}(x,Y_{x}(t_0,y),t).(D_xY_{x}(t_0,y).h).
\end{align*}

Then we see that we have to analyse
$$D_x\hat{F}(x,y,t).h-D_x\hat{F}(x,z,t).h$$ and $$D_y\hat{F}(x,y,t).$$

\begin{itemize}
\item For any $y,z\in H$, we have
\begin{align*}
|\hat{F}(x,y,t)-\hat{F}(x,z,t)|&=|\E F(x,Y_{x}(t,y))-\E F(x,Y_x(t,z))|\\
&\leq C\E|Y_{t}^{x}(y)-Y_{t}^{x}(z)|\\
&\leq Ce^{-ct}|y-z|,
\end{align*}

and we deduce that $|D_y\hat{F}(x,y,t).k|\leq Ce^{-ct}|k|$.
\item Moreover we know that $U_{t}^{x,h}(y)=D_xY_{x}(t,y).h$ is solution of
\begin{gather*}
dU_{t}^{x,h}(y)=\left(BU_{t}^{x,h}(y)+D_xG(x,Y_{x}(t,y)).h+D_yG(x,Y_x(t,y)).U_{t}^{x,h}(y)\right)dt\\
U_{0}^{x,h}(y)=0.
\end{gather*}

We deduce the following property: $|U_{t}^{x,h}(y)|\leq C|h|$ a.s.

As a consequence $|\E D_y\hat{F}(x,Y_x(t_0,y),t).(D_xY_x(t_0,y).h)|\leq Ce^{-ct}|h|$.

\item Now we take $x,y,z,t,h$, and we compute
\begin{align*}
D_x\hat{F}(x,y,t).h-&D_x\hat{F}(x,z,t).h=\E\left(D_xF(x,Y_x(t,y)).h-D_xF(x,Y_x(t,z)).h\right)\\
&+\E\left(D_yF(x,Y_x(t,y)).U_{t}^{x,h}(y)-D_yF(x,Y_x(t,z)).U_{t}^{x,h}(z)\right)\\
&=\E\left(D_xF(x,Y_x(t,y)).h-D_xF(x,Y_x(t,z)).h\right)\\
&+\E\left([D_yF(x,Y_x(t,y))-D_yF(x,Y_x(t,z))].U_{t}^{x,h}(y)\right)\\
&+\E\left(D_yF(x,Y_x(t,z)).(U_{t}^{x,h}(y)-U_{t}^{x,h}(z))\right)
\end{align*}
First, we have
\begin{align*}
|\E\left(D_xF(x,Y_x(t,y)).h-D_xF(x,Y_x(t,z)).h\right)|&\leq\E|D_xF(x,Y_x(t,y)).h-D_xF(x,Y_x(t,z)).h|\\
&\leq C|h|_{H}\E|Y_x(t,y)-Y_x(t,z)|_{(-B)^\eta}\\
&\leq C|h|_{H}e^{-ct}(1+\frac{1}{t^\eta})|y-z|_{H},
\end{align*}
using the regularity assumptions on $F$ (see \eqref{hypF}), and using the following estimate:
$$\E|Y_x(t,y)-Y_x(t,z)|_{(-B)^\eta}\leq Ce^{-ct}(1+\frac{1}{t^\eta})|y-z|_{H},$$
for some $c>0$.

Second,
\begin{align*}
|\E([D_yF(x,Y_x(t,y))-&D_yF(x,Y_x(t,z))].U_{t}^{x,h}(y))|\\
&\leq \E|[D_yF(x,Y_x(t,y))-D_yF(x,Y_x(t,z))].U_{t}^{x,h}(y)|\\
&\leq C\E|U_{t}^{x,h}(y)|_{H}|Y_x(t,y)-Y_x(t,z)|_{(-B)^\eta}\\
&\leq Ce^{-ct}(1+\frac{1}{t^\eta})|h|_{H}|y-z|_{H}.
\end{align*}

Third,
\begin{align*}
|\E\left(D_yF(x,Y_x(t,z)).(U_{t}^{x,h}(y)-U_{t}^{x,h}(z))\right)|&\leq \E|D_yF(x,Y_x(t,z)).(U_{t}^{x,h}(y)-U_{t}^{x,h}(z))|\\
&\leq C\E|U_{t}^{x,h}(y)-U_{t}^{x,h}(z)|_H;
\end{align*}

It remains to look at $|U_{t}^{x,h}(y)-U_{t}^{x,h}(z)|_H$; we indeed have
$$|U_{t}^{x,h}(y)-U_{t}^{x,h}(z)|^2\leq C|h|_{H}^{2}|y-z|_{H}^{2}e^{-c_0t},$$
where $c_0>0$.

We use these inequalities with $z:=Y_{x}(t_0,y)$; recalling that for any $t_0$
$$\E|Y_x(t_0,y)|\leq C(1+|y|),$$ we get
$$|D_x\tilde{F}_{t_0}(x,y,t).h|\leq  C(1+|y|)e^{-ct}(1+\frac{1}{t^\eta})|h|.$$

\end{itemize}
\qed

\end{document}